\documentclass{article}

\usepackage{microtype}
\usepackage{graphicx}
\usepackage{subfigure}
\usepackage{booktabs} 

\usepackage{hyperref}

\usepackage[accepted]{icml2021}

\icmltitlerunning{A Hybrid Variance-Reduced Method for Decentralized Stochastic Non-Convex Optimization}

\usepackage{bbm}
\usepackage{amsmath}   
\usepackage{amsthm}
\usepackage{appendix}
\usepackage{tabularx,amsmath,amssymb,graphicx,paralist,subfigure,pdfpages,paralist,multirow,comment}
\newtheorem{lem}{Lemma} 
\newtheorem{theorem}{Theorem}

\newtheorem{cor}{Corollary}
\newtheorem{rmk}{Remark}
\newtheorem{assump}{Assumption}

\usepackage{makecell}
\usepackage{verbatim}
\usepackage{pifont}
\newcommand{\cmark}{\text{\ding{51}}}
\newcommand{\xmark}{\text{\ding{55}}}

\def\mbb{\mathbb}
\def\mb{\mathbf}
\def\mc{\mathcal}
\def\tsum{\textstyle{\sum}}
\def\wh{\widehat}
\def\wt{\widetilde}
\def\ol{\overline}
\def\ul{\underline}

\def\bds{\boldsymbol}

\def\E{\mbb{E}}

\def\F{\mc{F}}
\def\R{\mathbb{R}}
\def\P{\mathbb{P}}
\def\HS{\texttt{GT-HSGD}}
\def\n{\nonumber}
\def\ra{\rightarrow}
\def\x{\mb x}
\def\y{\mb y}
\def\z{\mb z}
\def\s{\mb s}
\def\X{\bds\xi}
\def\W{\mb{W}}
\def\v{\mb{v}}
\def\g{\mb{g}}
\def\I{\mb{I}}
\def\J{\mb{J}}
\def\a{\alpha}
\def\f{\mb f}
\def\H{\mc H}
\def\SFO{\texttt{SFO}}
\def\DSGD{\texttt{DSGD}}

\renewcommand{\arraystretch}{1.5}

\begin{document}

\twocolumn[

           
\icmltitle{A Hybrid Variance-Reduced Method for \\ Decentralized Stochastic Non-Convex Optimization}



\icmlsetsymbol{equal}{*}

\begin{icmlauthorlist}
\icmlauthor{Ran Xin}{cmu}
\icmlauthor{Usman A. Khan}{tufts}
\icmlauthor{Soummya Kar}{cmu}
\end{icmlauthorlist}


\icmlaffiliation{cmu}{Department of Electrical and Computer Engineering, Carnegie Mellon University, Pittsburgh, PA, USA}

\icmlaffiliation{tufts}{Department of Electrical and Computer Engineering, Tufts University, Medford, MA, USA}


\icmlcorrespondingauthor{Ran Xin}{ranx@andrew.cmu.edu}

\icmlkeywords{Machine Learning, ICML}

\vskip 0.3in
]



\printAffiliationsAndNotice{}  


\begin{abstract}
This paper considers decentralized stochastic optimization over a network of~$n$ nodes, where each node possesses a smooth non-convex local cost function and the goal of the networked nodes is to find an~$\epsilon$-accurate first-order stationary point of the sum of the local costs. We focus on an online setting, where each node accesses its local cost only by means of a stochastic first-order oracle that returns a noisy version of the exact gradient. In this context, we propose a novel single-loop decentralized hybrid variance-reduced stochastic gradient method, called \texttt{GT-HSGD}, that outperforms the existing approaches in terms of both the oracle complexity and practical implementation. The \texttt{GT-HSGD} algorithm implements specialized local hybrid stochastic gradient estimators that are fused over the network to track the global gradient. Remarkably, \texttt{GT-HSGD} achieves a network topology-independent oracle complexity of~$O(n^{-1}\epsilon^{-3})$ when the required error tolerance~$\epsilon$ is small enough, leading to a linear speedup with respect to the centralized optimal online variance-reduced approaches that operate on a single node. Numerical experiments are provided to illustrate our main technical results.
\end{abstract}

\section{Introduction}
We consider~$n$ nodes, such as machines or edge devices, communicating over a decentralized network described by a directed graph~${\mc{G} = (\mc{V},\mc{E})}$, where~${\mc{V} = \{1,\cdots,n\}}$ is the set of node indices and~${\mc{E}}\subseteq\mc{V}\times\mc{V}$ is the collection of ordered pairs~${(i,j)}$,~${i,j \in \mc{V}}$, such that node~$j$ sends information to node~$i$. Each node~$i$ possesses a private local cost function~${f_i:\R^p \ra\R}$ and the goal of the networked nodes is to solve, via local computation and communication, the following optimization problem: $$\min_{\x\in\R^p}F(\x) = \frac{1}{n}\sum_{i=1}^n f_i(\x).$$
This canonical formulation is known as decentralized optimization~\cite{DGD_tsitsiklis,DGD_nedich,DGD_Kar,DSGD1_Chen} that has emerged as a promising framework for large-scale data science and machine learning problems~\cite{DSGD_NIPS,SGP_ICML}. Decentralized optimization is essential in scenarios where data is geographically distributed and/or centralized data processing is infeasible due to communication and computation overhead or data privacy concerns. 
In this paper, we focus on an \emph{online and non-convex} setting. In particular, we assume that each local cost~$f_i$ is \emph{non-convex} and each node~$i$ only accesses~$f_i$ by querying a local \emph{stochastic first-order oracle (\SFO)}~\cite{SGD_Nemirovski} that returns a stochastic gradient, i.e., a noisy version of the exact gradient, at the queried point. As a concrete example of practical interest, the~\SFO~mechanism applies to many online learning and expected risk minimization problems where the noise in~\SFO~lies in the uncertainty of sampling from the underlying streaming data received at each node~\cite{DGD_Kar,DSGD1_Chen}. We are interested in the oracle complexity, i.e., the total number of queries to~\SFO~required at each node, to find an~$\epsilon$-accurate first-order stationary point~$\x^*$ of the global cost~$F$ such that~$\E[\|\nabla F(\x^*)\|]\leq\epsilon$.

\subsection{Related Work}
We now briefly review the literature of decentralized non-convex optimization with \texttt{SFO}, which has been widely studied recently. Perhaps the most well-known approach is the decentralized stochastic gradient descent (\DSGD) and its variants~\cite{DSGD1_Chen,DGD_Kar,DSGD_vlaski_2,DSGD_NIPS,QSGP}, which combine average consensus and a local stochastic gradient step. Although being simple and effective, \DSGD~is known to have difficulties in handling heterogeneous data~\cite{SPM_Xin}. Recent works~\cite{D2,GNSD,improved_DSGT_Xin,PD_SGD} achieve robustness to heterogeneous environments by leveraging certain decentralized bias-correction techniques such as \texttt{EXTRA} (type)~\cite{EXTRA,SED,EXTRA_revisit}, gradient tracking~\cite{NEXT_scutari,GT_CDC,MP_Pu,PIEEE_Xin,DIGing,harnessing,add-opt}, and primal-dual principles~\cite{GT_jakovetic,NIDS,DGT_NIPS,PD_Xu}. Built on top of these bias-correction techniques, very recent works~\cite{D_Get} and~\cite{D-SPIDER-SFO} propose \texttt{D-GET} and \texttt{D-SPIDER-SFO} respectively that further incorporate online \texttt{SARAH}/\texttt{SPIDER}-type variance reduction schemes~\cite{spider,spiderboost,sarah_ncvx} to achieve lower oracle complexities, when the~$\SFO$ satisfies a mean-squared smoothness property.
Finally, we note that the family of decentralized variance reduced methods has been significantly enriched recently, see, for instance,~\cite{DSA,DAVRG,SPM_Xin,AccDVR,Network-DANE,DPD_VR,GT-SARAH,GTSAGA_NCVX,DVR_LHQ}; however, these approaches are explicitly designed for empirical minimization where each local cost~$f_i$ is decomposed as a finite-sum of component functions, i.e.,~${f_i = \frac{1}{m}\sum_{r=1}^mf_{i,r}}$; it is therefore unclear whether these algorithms can be adapted to the online~\SFO~setting, which is the focus of this paper.

\subsection{Our Contributions}
In this paper, we propose~\HS, a novel online variance reduced method for decentralized non-convex optimization with stochastic first-order oracles~(\SFO). To achieve fast and robust performance, the~\HS~algorithm is built upon global gradient tracking~\cite{NEXT_scutari,GT_CDC} and a local hybrid stochastic gradient estimator~\cite{HSARAH,HSARAH_0,mSARAH} that can be considered as a convex combination of the vanilla stochastic gradient returned by the~\SFO~and a \texttt{SARAH}-type variance-reduced stochastic gradient~\cite{SARAH}. In the following, we emphasize the key advantages of~\HS~compared with the existing decentralized online (variance-reduced) approaches, from both theoretical and practical aspects. 

\textbf{Improved oracle complexity}.
A comparison of the oracle complexity of~\HS~with related algorithms is provided in Table~\ref{comp_oc}, from which we have the following important observations.
First of all, the oracle complexity of \texttt{GT-HSGD} is lower than that of \texttt{DSGD}, \texttt{D2}, \texttt{GT-DSGD} and \texttt{D-PD-SGD}, which are decentralized online algorithms without variance reduction; however,~\HS~imposes on the \texttt{SFO} an additional mean-squared smoothness (MSS) assumption that is required by all online variance-reduced techniques in the literature~\cite{lowerbound_sgd,spider,spiderboost,sarah_ncvx,HSARAH,HSARAH_0,mSARAH,D_Get,D-SPIDER-SFO,SNVRG}.
Secondly, \texttt{GT-HSGD} further achieves a lower oracle complexity than the existing decentralized online variance-reduced methods \texttt{D-GET}~\cite{D_Get} and \texttt{D-SPIDER-SFO}~\cite{D-SPIDER-SFO}, especially in a regime where the required error tolerance~$\epsilon$ and the network spectral gap~$(1-\lambda)$ are relatively small.\footnote{A small network spectral gap~$(1-\lambda)$ implies that the connectivity of the network is weak.}.
Moreover, when~$\epsilon$ is small enough such that $\epsilon\lesssim\min\big\{\lambda^{-4}(1-\lambda)^3n^{-1},\lambda^{-1}(1-\lambda)^{1.5}n^{-1}\big\}$, it can be verified that the oracle complexity of~\HS~reduces to~$O(n^{-1}\epsilon^{-3})$, independent of the network topology, and~\HS~achieves a linear speedup, in terms of the scaling with the network size~$n$, compared with the centralized optimal online variance-reduced approaches that operate on a single node~\cite{spider,spiderboost,sarah_ncvx,HSARAH,HSARAH_0,SNVRG}; see Section~\ref{sec_main_results} for a detailed discussion. In sharp contrast, the speedup of \texttt{D-GET}~\cite{D_Get} and \texttt{D-SPIDER-SFO}~\cite{D-SPIDER-SFO} is not clear compared with the aforementioned centralized optimal methods even if the network is fully connected, i.e.,~$\lambda = 0$.     
   
\textbf{More practical implementation}. Both \texttt{D-GET}~\cite{D_Get} and \texttt{D-SPIDER-SFO}~\cite{D-SPIDER-SFO} are double-loop algorithms that require very large minibatch sizes. In particular, during each inner loop they execute a fixed number of minibatch stochastic gradient type iterations with~$O(\epsilon^{-1})$ oracle queries per update per node, while at every outer loop they obtain a stochastic gradient with mega minibatch size by~$O(\epsilon^{-2})$ oracle queries at each node. Clearly,
querying the oracles exceedingly, i.e., obtaining a large amount of samples, at each node and every iteration in online steaming data scenarios substantially jeopardizes the actual wall-clock time. This is because the next iteration cannot be performed until all nodes complete the sampling process. Moreover, the double-loop implementation may incur periodic network synchronizations. These issues are especially significant when the working environments of the nodes are heterogeneous. Conversely, the proposed~\HS~is a single-loop algorithm with~$O(1)$ oracle queries per update and only requires a large minibatch size with~$O(\epsilon^{-1})$ oracle queries once in the \emph{initialization phase}, i.e., before the update recursion is executed; see Algorithm~\ref{GTH} and Corollary~\ref{main} for details. 

\renewcommand{\arraystretch}{2}
\begin{table*}[!ht]
\footnotesize
\caption{A comparison of the oracle complexity of decentralized online stochastic gradient methods. The oracle complexity is in terms of the total number of queries to~\SFO~required \emph{at each node} to obtain an~$\epsilon$-accurate stationary point~$\x^*$ of the global cost~$F$ such that~$\E[\|\nabla F(\x^*)\|]\leq\epsilon$. In the table,~$n$ is the number of the nodes and~$(1 - \lambda)\in(0,1]$ is the spectral gap of the weight matrix associated with the network. We note that the complexity of \texttt{D2} and \texttt{D-SPIDER-SFO} also depends on the smallest eigenvalue~$\lambda_n$ of the weight matrix; however, since~$\lambda_n$ is less sensitive to the network topology, we omit the dependence of~$\lambda_n$ in the table for conciseness. 
The MSS column indicates whether the algorithm in question requires the mean-squared smoothness assumption on the \texttt{SFO}. Finally, we emphasize that \texttt{DSGD} requires bounded heterogeneity such that~$\sup_{\x}\frac{1}{n}\sum_{i=1}^n\|\nabla f_i(\x) - \nabla F(\x)\|^2\leq\zeta^2$, for some~$\zeta\in\R^{+}$, while other algorithms in the table do not need this assumption.}
\label{comp_oc}
\begin{center}
\begin{tabular}{|c|l|c|c|r|}
\hline
\textbf{Algorithm} & \textbf{Oracle Complexity} & \textbf{MSS} & \textbf{Remarks}\\ \hline
\texttt{DSGD}~\cite{DSGD_NIPS} & $O\left(\max\left\{\dfrac{1}{n\epsilon^{4}}, \dfrac{\lambda^2n}{(1-\lambda)^2\epsilon^2}\right\} \right)$ & \xmark & bounded heterogeneity  \\ \hline
\texttt{D2}~\cite{D2} & $O\left(\max\left\{\dfrac{1}{n\epsilon^{4}} , \dfrac{n}{(1-\lambda)^b\epsilon^2}\right\}\right)$ & \xmark & \makecell{$b\in\R^{+}$ is not explicitly  \\ shown in~\cite{D2}} \\ \hline
\texttt{GT-DSGD}~\cite{improved_DSGT_Xin} & $O\left(\max\left\{\dfrac{1}{n\epsilon^{4}}, \dfrac{\lambda^2n}{(1-\lambda)^3\epsilon^2}\right\}\right)$ & \xmark &  \\ \hline
\texttt{D-PD-SGD}~\cite{PD_SGD} & $O\left(\max\left\{\dfrac{1}{n\epsilon^{4}}, \dfrac{n}{(1-\lambda)^{c}\epsilon^2}\right\}\right)$ & \xmark & \makecell{$c\in\R^{+}$ is not explicitly \\  shown in~\cite{PD_SGD}}  \\ \hline
\texttt{D-GET}~\cite{D_Get} &  $O\left(\dfrac{1}{(1-\lambda)^d\epsilon^{3}}\right)$ & \cmark & \makecell{$d\in\R^{+}$ is not explicitly \\ shown in~\cite{D_Get}} \\ \hline 
\makecell{\texttt{D-SPIDER-SFO}\\\cite{D-SPIDER-SFO}} &  $O\left(\dfrac{1}{(1-\lambda)^h\epsilon^{3}}\right)$ & \cmark & \makecell{$h\in\R^{+}$ is not explicitly \\  shown in~\cite{D-SPIDER-SFO}}\\ \hline 
\makecell{\texttt{GT-HSGD} \\ (\textbf{this work})} &  $O\left(\max\left\{\dfrac{1}{n\epsilon^{3}},
    \dfrac{\lambda^4}{(1-\lambda)^3\epsilon^2},  \dfrac{\lambda^{1.5}n^{0.5}}{(1-\lambda)^{2.25}\epsilon^{1.5}}\right\}\right)$ & \cmark &  \\ \hline
\end{tabular}
\end{center}
\end{table*}

\subsection{Roadmap and Notations}
The rest of the paper is organized as follows.
In Section~\ref{sec_problem_alg}, we state the problem formulation and develop the proposed \texttt{GT-HSGD} algorithm. Section~\ref{sec_main_results} presents the main convergence results of~\HS~and their implications. Section~\ref{sec_conv_analysis} outlines the convergence analysis of \texttt{GT-HSGD}, while the detailed proofs are provided in the Appendix.  Section~\ref{sec_num_exp} provides numerical experiments to illustrate our theoretical claims. Section~\ref{sec_conc} concludes the paper. 

We adopt the following notations throughout the paper.
We use lowercase bold letters to denote vectors and uppercase bold letters to denote matrices. The ceiling function is denoted as~$\lceil\cdot\rceil$.
The matrix~$\mb{I}_d$ represents the~$d\times d$ identity; $\mb{1}_d$ and $\mb{0}_d$ are the~$d$-dimensional column vectors of all ones and zeros, respectively. We denote~$[\x]_i$ as the~$i$-th entry of a vector~$\x$.
The Kronecker product of two matrices~$\mb{A}$ and~$\mb{B}$ is denoted by~$\mb{A}\otimes \mb{B}$. We use~$\|\cdot\|$ to denote the Euclidean norm of a vector or the spectral norm of a matrix. We use $\sigma(\cdot)$ to denote the~$\sigma$-algebra generated by the sets and/or random vectors in its argument.

\section{Problem Setup and \HS}\label{sec_problem_alg}
In this section, we introduce the mathematical model of the stochastic first-order oracle (\texttt{SFO}) at each node and the communication network. Based on these formulations, we develop the proposed~\HS~algorithm.

\subsection{Optimization and Network Model} 
We work with a rich enough probability space~$\{\Omega,\P,\F\}$. We consider decentralized recursive algorithms of interest that generate a sequence of estimates~$\{\x_t^i\}_{t\geq0}$ of the first-order stationary points of~$F$ at each node~$i$, where~$\x_0^i$ is assumed constant.  
At each iteration~$t$, each node~$i$ observes a random vector~$\X_{t}^i$ in~$\R^q$, which, for instance, may be considered as noise or as an online data sample. 
We then introduce the natural filtration (an increasing family of sub-$\sigma$-algebras of~$\F$) induced by these random vectors observed sequentially by the networked nodes:
\begin{align}\label{Ft}
\F_0 :=&~\{\Omega,\phi\}, \n\\
\F_t :=&~\sigma\left(\left\{\X_{0}^i,\X_{1}^i,\cdots,\X_{t-1}^i: i\in\mc{V}\right\}\right), \qquad\forall t\geq1,
\end{align}
where~$\phi$ is the empty set. We are now ready to define the \SFO~mechanism in the following. At each iteration~$t$,
each node~$i$, given an input random vector~$\x\in\R^p$ that is~$\F_t$-measurable, is able to query the local \SFO~to obtain a stochastic gradient of the form~$\g_i(\x,\X_{t}^i)$, where~$\g_i:\R^p\times\R^q\rightarrow\R^p$ is a Borel measurable function. 
We assume that the \texttt{SFO} satisfies the following four properties. 
\begin{assump}[\textbf{Oracle}]\label{oracle_asp}
\normalfont
For any~$\F_t$-measurable random vectors~$\x,\y\in\R^p$, we have the following:~$\forall i\in\mc{V}$,~$\forall t\geq0$,
\begin{itemize}
\item~$\E\left[\g_i(\x,\X_{t}^i)|\F_t\right] = \nabla f_i(\x)$;
\item~$\E\left[\|\g_i(\x,\X_{t}^i)-\nabla f_i(\x)\|^2\right] \leq\nu_i^2,$~$\ol{\nu}^2 := \frac{1}{n}\sum_{i=1}^n\nu_i^2$;
\item the family~$\left\{\X_{t}^i:\forall t\geq0,i\in\mc{V}\right\}$ of random vectors is independent;
\item~$\E\left[\|\g_i(\x,\X_{t}^i) - \g_i(\y,\X_{t}^i)\|^2\right]\leq L^2\E\left[\|\x-\y\|^2\right]$.
\end{itemize}
\end{assump}
The first three properties above are standard and commonly used to establish the convergence of decentralized stochastic gradient methods. They however do not explicitly impose any structures on the stochastic gradient mapping~$\g_i$ other than the measurability. On the other hand, the last property, the mean-squared smoothness, roughly speaking, requires that~$\g_i$ is~$L$-smooth on average with respect to the input arguments~$\x$ and~$\y$. As a simple example, Assumption~\ref{oracle_asp} holds if~${f_i(\x) = \frac{1}{2}\x^{\top}\mb{Q}_i\x}$ and~${\g_i(\x,\X_i) = \mb{Q}_i\x + \X_i}$, where~$\mb{Q}_i$ is a constant matrix and~$\X_i$ has zero mean and finite second moment.
We further note that the mean-squared smoothness of each~$\g_i$ implies, by Jensen's inequality, that each~$f_i$ is~$L$-smooth, i.e.,~$\|\nabla f_i(\x) - \nabla f_i(\y)\|\leq L\|\x-\y\|$, and consequently the global function~$F$ is also~$L$-smooth. 

In addition, we make the following assumptions on~$F$ and the communication network~$\mc{G}$.
\begin{assump}[\textbf{Global Function}]
\normalfont
$F$ is bounded below, i.e.,~$F^* := \inf_{\x\in\R^p}F(\x)>-\infty$.
\end{assump}
\begin{assump}[\textbf{Communication Network}]\label{net_asp}
\normalfont
The directed network~$\mc{G}$ admits a primitive and doubly-stochastic weight matrix~$\ul{\W} = \{\ul{w}_{ij}\}\in\R^{n\times n}$. Hence,~$\ul{\W}\mb{1}_n = \ul{\W}^\top\mb{1}_n = \mb{1}_n$ and~$\lambda := \|\ul{\W} - \frac{1}{n}\mb{1}_n\mb{1}_n^\top\| \in[0,1)$.
\end{assump} 

The weight matrix~$\ul{\W}$ that satisfies Assumption~\ref{net_asp} may be designed for strongly-connected weight-balanced directed graphs (and thus for arbitrary connected undirected graphs). For example, the family of directed exponential graphs is weight-balanced and plays a key role in decentralized training~\cite{SGP_ICML}.
We note that~$\lambda$ is known as the second largest singular value of~$\ul{\W}$ and measures the algebraic connectivity of the graph, i.e., a smaller value of~$\lambda$ roughly means a better connectivity.
We note that several existing approaches require strictly stronger assumptions on~$\ul{\W}$. For instance,~\texttt{D2}~\cite{D2} and~\texttt{D-PD-SGD}~\cite{PD_SGD} require~$\ul{\W}$ to be symmetric and hence are restricted to undirected networks.

\subsection{Algorithm Development and Description}
We now describe the proposed~\HS~algorithm and provide an intuitive construction. Recall that~$\x_t^i$ is the estimate of an stationary point of the global cost~$F$ at node~$i$ and iteration~$t$. Let~$\g_i(\x_t^i,\X_{t}^i)$ and~$\g_i(\x_{t-1}^i,\X_{t}^i)$  be the corresponding stochastic gradients returned by the local \SFO~queried at~$\x_t^i$ and~$\x_{t-1}^i$ respectively.  Motivated by the strong performance of recently introduced decentralized methods that combine gradient tracking and various variance reduction schemes for finite-sum problems~\cite{GT-SARAH,GTSAGA_NCVX,Network-DANE,D_Get}, we seek similar variance reduction for decentralized online problems with~\SFO. 
In particular, we focus on the following \textit{local} hybrid variance reduced stochastic gradient estimator~$\mb{v}_t^i$ introduced in~\cite{HSARAH,HSARAH_0,mSARAH} for centralized online problems:~$\forall t\geq1$,
\begin{align}\label{vi}
\mb{v}_{t}^i 
= \g_i(\x_t^i,\X_{t}^i) + (1-\beta)\big(\mb{v}_{t-1}^i - \g_i(\x_{t-1}^i,\X_{t}^i)\big),
\end{align}
for some applicable weight parameter~$\beta\in[0,1]$. This local gradient estimator~$\mb{v}_t^i$ is fused, via a gradient tracking mechanism~\cite{NEXT_scutari,GT_CDC}, over the network to update the global gradient tracker~$\mb y_t^i$, which is subsequently used as the descent direction in the~$\mb x_t^i$-update. The complete description of~\HS~is provided in Algorithm~\ref{GTH}. We note that the update~\eqref{vi} of~$\mb{v}_t^i$ may be equivalently written as
\begin{align*}
\mb{v}_{t}^i 
=&~~\beta\cdot\!\!\!\!\!\!\!\!\!\!\underbrace{\g_i(\x_t^i,\X_{t}^i)}_{\texttt{Stochastic gradient}} \n\\
&~+(1-\beta)\cdot\underbrace{\left(\g_i(\x_t^i,\X_{t}^i) - \g_i(\x_{t-1}^i,\X_{t}^i) + \mb{v}_{t-1}^i\right)}_{\texttt{SARAH}},
\end{align*}
which is a convex combination of the local vanilla stochastic gradient returned by the \SFO~and a \texttt{SARAH}-type~\cite{SARAH,spider,spiderboost} gradient estimator. 
This discussion leads to the fact that~\HS~reduces to \texttt{GT-DSGD}~\cite{MP_Pu,improved_DSGT_Xin,GNSD} when~$\beta=1$, and becomes the inner loop of~\texttt{GT-SARAH}~\cite{GT-SARAH} when~$\beta=0$. However, our convergence analysis shows that \texttt{GT-HSGD} achieves its best oracle complexity and outperforms the existing decentralized online variance-reduced approaches~\cite{D_Get,D-SPIDER-SFO} with a weight parameter~${\beta\in(0,1)}$.
It is then clear that neither \texttt{GT-DSGD} nor the inner loop of \texttt{GT-SARAH}, on their own, are able to outperform the proposed approach, making~\HS~a non-trivial algorithmic design for this problem class.

\begin{algorithm}[tbph]
\caption{\HS~at each node~$i$}
\label{GTH}
\begin{algorithmic}[1]
\REQUIRE{$\mb{x}_0^i = \ol{\mb{x}}_0$;~$\alpha$;~$\beta$;~$b_0$;~$\mb{y}_0^i = \mb{0}_p$;~$\mb{v}_{-1}^i = \mb{0}_p$;~$T$.}
\vspace{0.1cm}
\STATE{Sample~$\{\X_{0,r}^{i}\}_{r=1}^{b_0}$ and~$\mb{v}_{0}^i = \frac{1}{b_0}\sum_{r=1}^{b_0}\g_i(\x_0^i,\X_{0,r}^i);$}
\vspace{0.1cm}
\STATE{$\mb{y}_{1}^{i} = \sum_{j=1}^{n}\ul{w}_{ij}\big(\mb{y}_{0}^j + \mb{v}_0^{j} - \mb{v}_{-1}^{j}\big);$}
\vspace{0.1cm}
\STATE{$\mb{x}_{1}^i = \sum_{j=1}^{n}\ul{w}_{ij}\big(\mb{x}_{0}^j - \alpha\mb{y}_{1}^j\big);$}
\vspace{0.1cm}
\FOR{{$t= 1,2,\cdots,T-1$}}
\vspace{0.1cm}
\STATE{Sample~$\X_{t}^i$;}
\vspace{0.1cm}
\STATE{$\mb{v}_{t}^i 
= \g_i(\x_t^i,\X_{t}^i) + (1-\beta)\big(\mb{v}_{t-1}^i - \g_i(\x_{t-1}^i,\X_{t}^i)\big).$}
\vspace{0.1cm}
\STATE{$\mb{y}_{t+1}^{i} = \tsum_{j=1}^{n}\ul{w}_{ij}\big(\mb{y}_{t}^j + \mb{v}_t^{j} - \mb{v}_{t-1}^{j}\big);$}
\vspace{0.1cm}
\STATE{$\mb{x}_{t+1}^i = \tsum_{j=1}^{n}\ul{w}_{ij}\big(\mb{x}_{t}^j - \alpha\mb{y}_{t+1}^j\big);$}
\vspace{0.1cm}
\ENDFOR
\OUTPUT{$\wt{\x}_{T}$ selected uniformly at random from~$\{\x_t^i\}_{0\leq t\leq T}^{i\in\mc{V}}$.}
\end{algorithmic}
\end{algorithm}

\begin{rmk}
Clearly, each~$\v_t^i$ is a conditionally biased estimator of~$\nabla f_i(\x_t^i)$, i.e.,~$\E[\v_t^i|\F_t] \neq \nabla f_i(\x_t^i)$ in general. 
However, it can be shown that~$\E[\v_t^i] = \E[\nabla f_i(\x_t^i)]$, meaning that~$\v_t^i$ serves as a surrogate for the underlying exact gradient in the sense of total expectation.
\end{rmk}

\section{Main Results}\label{sec_main_results}
In this section, we present the main convergence results of~\HS~in this paper and discuss their salient features. The formal convergence analysis is deferred to Section~\ref{sec_conv_analysis}.

\begin{theorem}\label{main0}
If the weight parameter~$\beta = \frac{48L^2\a^2}{n}$ and the step-size~$\a$ is chosen as
\begin{align*}
0<\a<\min\bigg\{\frac{(1-\lambda^2)^2}{90\lambda^2},\frac{\sqrt{n(1-\lambda)}}{26\lambda},\frac{1}{4\sqrt{3}}\bigg\}\frac{1}{L},
\end{align*}
then the output~$\wt{\x}_T$ of \texttt{GT-HSGD} satisfies:~$\forall T\geq2$,
\begin{align*}
\E\big[\|&\nabla F(\wt{\x}_T)\|^2\big] 
\leq 
\frac{4(F(\ol{\x}_0) - F^*)}{\a T}
+ \frac{8\beta\ol{\nu}^2}{n}
+ \frac{4\ol{\nu}^2}{\beta b_0nT} \n\\
&+ \frac{64\lambda^4\|\nabla\mb{f}\big(\x_{0}\big)\|^2}{(1-\lambda^2)^3nT}
+ \frac{96\lambda^2\ol{\nu}^2}{(1-\lambda^2)^{3}b_0T}
+\frac{256\lambda^2\beta^2\ol{\nu}^2}{(1-\lambda^2)^3},
\end{align*}
where~${\|\nabla\mb{f}\big(\x_{0}\big)\|^2 = \sum_{i=1}^n\|\nabla f_i\big(\ol{\x}_{0}\big)\|^2}$
\end{theorem}

\vspace{0.3cm}
\begin{rmk}
Theorem~\ref{main0} holds for~\HS~with arbitrary initial minibatch size~$b_0\geq1$.
\end{rmk}

Theorem~\ref{main0} establishes a non-asymptotic bound, with no hidden constants, on the mean-squared stationary gap of \texttt{GT-HSGD} over any finite time horizon~$T$. 

\textbf{Transient and steady-state performance over infinite time horizon.} If~$\a$ and~$\beta$ are chosen according to Theorem~\ref{main0}, the mean-squared stationary gap~$\E\left[\|\nabla F(\wt{\x}_{T})\|^2\right]$ of \texttt{GT-HSGD} decays sublinearly at a rate of~$O(1/T)$ up to a steady-state error (SSE) such that
\begin{align}\label{SSE}
\limsup_{T\ra\infty}\E\left[\|\nabla F(\wt{\x}_{T})\|^2\right]
\leq 
\frac{8\beta\ol{\nu}^2}{n}
+ \frac{256\lambda^2\beta^2\ol{\nu}^2}{(1-\lambda^2)^3}.
\end{align}
In view of~\eqref{SSE}, the SSE of~\HS~is bounded by the sum of two terms: (i) the first term is in the order of~$O(\beta)$ and the division by~$n$ demonstrates the benefit of increasing the network size\footnote{Since~\HS~computes~$O(n)$ stochastic gradients in parallel per iteration across the nodes, the network size~$n$ can be interpreted as the minibatch size of~\HS.}; (ii) the second term 
is in the order of~$O(\beta^2)$ and
reveals the impact of the spectral gap~$(1-\lambda)$ of the network topology. 
Clearly, the SSE can be made arbitrarily small by choosing small enough~$\beta$ and~$\a$. Moreover, since the spectral gap~$(1-\lambda)$ only appears in a higher order term of~$\beta$ in~\eqref{SSE}, its impact reduces as~$\beta$ becomes smaller, i.e., as we require a smaller SSE.

The following corollary is concerned with the finite-time convergence rate of \texttt{GT-HSGD} with specific choices of the algorithmic parameters~$\a,\beta$, and~$b_0$.

\begin{cor}\label{main}
Setting~$\a = \frac{n^{2/3}}{8LT^{1/3}}$,~$\beta = \frac{3n^{1/3}}{4T^{2/3}}$, and~$b_0 = \lceil \frac{T^{1/3}}{n^{2/3}}\rceil$ in Theorem~\ref{main0}, we have:
\begin{align*}
\E\big[\|\nabla F(\wt{\x}_T)\|^2\big] 
\leq&~\frac{32L(F(\ol{\x}_0) - F^*)+ 12\ol{\nu}^2}{(nT)^{2/3}}       \n\\
+& \frac{64\lambda^4\|\nabla\mb{f}\big(\x_{0}\big)\|^2}{(1-\lambda^2)^3nT}
+ \frac{240\lambda^2n^{2/3}\ol{\nu}^2}{(1-\lambda^2)^{3}T^{4/3}},
\end{align*}
for all~$T> \max\big\{\frac{1424\lambda^6n^{2}}{(1-\lambda^2)^6},\frac{35\lambda^3 n^{0.5}}{(1-\lambda)^{1.5}}\big\}$. As a consequence, \texttt{GT-HSGD} achieves an~$\epsilon$-accurate stationary point~$\mb{x}^*$ of the global cost~$F$ such that~$\E[\|\nabla F(\mb{x}^*)\|]\leq\epsilon$ with $$\mc{H} = O\left(\max\left\{\H_{opt},\H_{net}\right\}\right)$$ 
iterations\footnote{The~$O(\cdot)$ notation here does not absorb any problem parameters, i.e., it only hides universal constants.}, where~$\H_{opt}$ and~$\H_{net}$ are given respectively by
\begin{align*}
&\H_{opt} = \frac{(L(F(\ol{\x}_0) - F^*)+ \ol{\nu}^2)^{1.5}}{n\epsilon^{3}},
\\
&\H_{net} =\max\bigg\{\frac{\lambda^4\|\nabla\mb{f}\big(\x_{0}\big)\|^2}{(1-\lambda^2)^3n\epsilon^2},  \frac{\lambda^{1.5}n^{0.5}\ol{\nu}^{1.5}}{(1-\lambda^2)^{2.25}\epsilon^{1.5}}\bigg\}.
\end{align*}
The resulting total number of oracle queries at each node is thus~$\lceil\mc{H} + \mc{H}^{1/3}n^{-2/3}\rceil$.
\end{cor}

\vspace{0.1cm}
\begin{rmk}
Since~$\mc{H}^{1/3}n^{-2/3}$ is much smaller than~$\H$, we treat the oracle complexity of \texttt{GT-HSGD} as~$\H$ for the ease of exposition in Table~\ref{comp_oc} and the following discussion.
\end{rmk}
An important implication of Corollary~\ref{main} is given in the following.

\noindent\textbf{A regime for network topology-independent oracle complexity and linear speedup.} According to Corollary~\ref{main0}, the oracle complexity of~\HS~at each node is bounded by the maximum of two terms: (i) the first term~$\H_{opt}$ is independent of the network topology and, more importantly, is~$n$ times smaller than the oracle complexity of the optimal centralized online variance-reduced methods that execute on a single node for this problem class~\cite{spider,spiderboost,sarah_ncvx,HSARAH_0,HSARAH}; (ii) the second term~$\H_{net}$ depends on the network spectral gap~$1-\lambda$ and is in the lower order of~$1/\epsilon$. These two observations lead to the interesting fact that the oracle complexity of~\HS~becomes independent of the network topology, i.e.,~$\H_{opt}$ dominates~$\H_{net}$, if the required error tolerance~$\epsilon$ is small enough such that\footnote{This boundary condition follows from basic algebraic manipulations.}~$\epsilon\lesssim\min\big\{\lambda^{-4}(1-\lambda)^3n^{-1},\lambda^{-1}(1-\lambda)^{1.5}n^{-1}\big\}$. In this regime,  \texttt{GT-HSGD} thus achieves a network topology-independent oracle complexity~$\H_{opt} = O(n^{-1}\epsilon^{-3})$, exhibiting a linear speed up compared with the aforementioned centralized optimal algorithms~\cite{spider,spiderboost,sarah_ncvx,HSARAH_0,HSARAH,SNVRG}, in the sense that the total number of oracle queries required to achieve an~$\epsilon$-accurate stationary point at each node is reduced by a factor of~$1/n$.

\begin{rmk}\label{regime}
The small error tolerance regime in the above discussion corresponds to a large number of oracle queries, which translates to the scenario where the required total number of iterations~$T$ is large. Note that a large~$T$ further implies that the step-size~$\a$ and~the weight parameter~$\beta$ are small; see the expression of~$\a$ and~$\beta$ in Corollary~\ref{main}. 
\end{rmk}

\section{Outline of the Convergence Analysis}\label{sec_conv_analysis}
In this section, we outlines the proof of Theorem~\ref{main0}, while the detailed proofs are provided in the Appendix.
We let Assumptions~\ref{oracle_asp}-\ref{net_asp} hold throughout the rest of the paper without explicitly stating them. For the ease of exposition, we write the~$\x_t$- and~$\y_t$-update of~\HS~in the following equivalent matrix form:~$\forall t\geq0$,
\begin{subequations}
\begin{align}
\y_{t+1} &= \W\left(\mb{y}_{t} + \mb{v}_{t} - \mb{v}_{t-1}\right), \label{y_up}\\ 
\x_{t+1} &= \W\left(\x_t - \alpha\mb{y}_{t+1}\right),     \label{x_up}
\end{align}
\end{subequations}
where~$\W := \ul{\W}\otimes \I_p$ and~$\x_t, \y_t, \mb{v}_t$ are square-integrable random vectors in~$\mathbb{R}^{np}$ that respectively concatenate the local estimates~$\{\x_t^i\}_{i=1}^n$ of a stationary point of~$F$, gradient trackers~$\{\y_t^i\}_{i=1}^n$, stochastic gradient estimators~$\{\mb{v}_t^i\}_{i=1}^n$. 
It is straightforward to verify that~$\x_t$ and~$\y_t$ are~$\F_t$-measurable while~$\mb{v}_t$ is~$\F_{t+1}$-measurable for all~$t\geq0$.
For convenience, we also denote~$$\nabla\mb{f}(\mb{x}_t) := \left[\nabla f_1(\mb{x}_t^1)^\top,\cdots, 
\nabla f_n(\mb{x}_t^n)^\top\right]^\top$$ and introduce the following quantities:
\begin{align*}
&\mb{J} := \left(\frac{1}{n}\mb{1}_n\mb{1}_n^\top\right)\otimes \I_p \\
&\ol{\mb x}_t := \frac{1}{n}(\mb{1}_n^\top\otimes \I_p)\mb{x}_t,
\quad
\ol{\mb y}_t := \frac{1}{n}(\mb{1}_n^\top\otimes \I_p)\mb{y}_t, \\
&\ol{\mb v}_t := \frac{1}{n}(\mb{1}_n^\top\otimes \I_p)\mb{v}_t,
\quad
\ol{\nabla\mb{f}}(\mb{x}_t) := \frac{1}{n}(\mb{1}_n^\top\otimes \I_p)\nabla\mb{f}(\mb{x}_t).
\end{align*}

In the following lemma, we enlist several well-known results in the context of gradient tracking-based algorithms for decentralized stochastic optimization, whose proofs may be found in~\cite{NEXT_scutari,harnessing,GTVR,MP_Pu}.
\begin{lem}\label{basic}
The following relationships hold.
\begin{enumerate}[(a)]
\item $\|\W\mb{x}-\J\mb{x}\|\leq\lambda\|\mb{x}-\J\mb{x}\|,\forall\mb{x}\in\mbb{R}^{np}$. \label{W}
\item $\ol{\mb y}_{t+1} = \ol{\mb v}_t ,\forall t\geq0$. \label{track}
\item $\left\|\ol{\nabla\mb{f}}(\mb{x}_t)-\nabla F(\ol{\mb{x}}_t)\right\|^2\leq\frac{L^2}{n}\left\|\mb{x}_t-\J\mb{x}_t\right\|^2,\forall t\geq0$. \label{Lbound}
\end{enumerate}
\end{lem}
We note that Lemma~\ref{basic}(\ref{W}) holds since~$\ul{\W}$ is primitive and doubly-stochastic, Lemma~\ref{basic}(\ref{track}) is a direct consequence of the gradient tracking update~\eqref{y_up} and Lemma~\ref{basic}(\ref{Lbound}) is due to the~$L$-smoothness of each~$f_i$.
By the estimate update of~\HS~described in~\eqref{x_up} and Lemma~\ref{basic}(\ref{track}), it is straightforward to obtain: 
\begin{align}\label{ave}
\ol{\x}_{t+1} 
= \ol{\x}_t - \a \ol{\y}_{t+1} 
= \ol{\x}_t - \a \ol{\mb{v}}_t, \qquad\forall t\geq0. 
\end{align}
Hence, the mean state~$\ol{\x}_t$ proceeds in the direction of the average of local stochastic gradient estimators~$\ol{\mb{v}}_t$. 
With the help of~\eqref{ave} and the~$L$-smoothness of~$F$ and each~$f_i$,
we establish the following descent inequality which sheds light on the overall convergence analysis. 
\begin{lem}\label{DS0}
If~$0<\alpha\leq\frac{1}{2L}$, then we have:~$\forall T\geq0$,
\begin{align*}
\sum_{t=0}^{T}\big\|&\nabla F(\ol{\x}_{t})\big\|^2
\leq 
\frac{2(F(\ol{\x}_{0}) - F^*)}{\a}
- \frac{1}{2}\sum_{t=0}^{T}\left\|\ol{\mb{v}}_{t}\right\|^2\\
&+ 2\sum_{t=0}^{T}\left\|\ol{\mb{v}}_{t}-\ol{\nabla\mb{f}}(\x_{t})\right\|^2 
+ \frac{2L^2}{n}\sum_{t=0}^{T}\left\|\mb{x}_{t}-\J\x_{t}\right\|^2.
\end{align*}
\end{lem}


In light of Lemma~\ref{DS0}, our approach to establishing the convergence of~\HS~is to seek the conditions on the algorithmic parameters of~\HS, i.e., the step-size~$\a$ and the weight parameter~$\beta$, such that
\begin{align}\label{crt}
&- \frac{1}{2T}\sum_{t=0}^{T}\E\left[\|\ol{\mb{v}}_{t}\|^2\right] 
+ \frac{2}{T}\sum_{t=0}^{T}\E\left[\|\ol{\mb{v}}_{t}-\ol{\nabla\mb{f}}(\x_{t})\|^2\right]  \nonumber\\
&+ \frac{2L^2}{nT}\sum_{t=0}^{T}\E\left[\|\mb{x}_{t}-\J\x_{t}\|^2\right] 
= O\left(\a,\beta,\frac{1}{b_0},\frac{1}{T}\right),
\end{align}
where~$O(\a,\beta,1/b_0,1/T)$ represents a nonnegative quantity which may be made arbitrarily small by choosing small enough $\a$ and~$\beta$ along with large enough~$T$ and~$b_0$. If~\eqref{crt} holds, then Lemma~\ref{DS0} reduces to
\begin{align*}
\frac{1}{T+1}\sum_{t=0}^{T}&\E\left[\|\nabla F(\ol{\x}_{t})\|^2\right] \nonumber\\
&\leq 
\frac{2(F(\ol{\x}_{0}) - F^*)}{\a T}
+ O\left(\a,\beta,\frac{1}{b_0},\frac{1}{T}\right),
\end{align*}
which leads to the convergence arguments of~\HS. For these purposes, we quantify~$\sum_{t=0}^{T}\E\left[\|\ol{\v}_{t}-\ol{\nabla\mb{f}}(\x_{t})\|^2\right]$ and~$\sum_{t=0}^{T}\E\left[\|\mb{x}_{t}-\J\x_{t}\|^2\right]$ next.

\subsection{Contraction Relationships}\label{Aux}
First of all, we establish upper bounds on the gradient variances~$\E\left[\|\ol{\mb{v}}_t-\ol{\nabla\f}(\x_t)\|^2\right]$ and~$\E\left[\|\mb{v}_t-\nabla\f(\x_t)\|^2\right]$ by exploiting the hybrid and recursive update of~$\mb{v}_t$.  
\begin{lem}\label{vrb}
The following inequalities hold:~$\forall t\geq1$,
\begin{align}
&\E\left[\|\ol{\mb{v}}_{t} - \ol{\nabla \f}(\x_t)\|^2\right] \n\\
\leq&~(1-\beta)^2\E\left[\|\ol{\mb{v}}_{t-1} - \ol{\nabla \f}(\x_{t-1})\|^2\right] \n\\
&+\frac{6L^2\a^2}{n}(1-\beta)^2\E\left[\|\ol{\mb{v}}_{t-1}\|^2\right] 
+\frac{2\beta^2\ol{\nu}^2}{n}\n\\
&+\frac{6L^2}{n^2}(1-\beta)^2\E\left[\|\x_t - \J\x_t\|^2 + \|\x_{t-1}- \J\x_{t-1}\|^2\right],
\label{vrb_a}
\end{align}
and,~$\forall t\geq1$,
\begin{align}
&\E\left[\|\mb{v}_t - \nabla\f(\x_t)\|^2\right]
\n\\
\leq&~(1-\beta)^2\E\left[\|\mb{v}_{t-1} - \nabla\f(\x_{t-1})\|^2\right]  \n\\
&+ 6nL^2\a^2(1-\beta)^2\E\left[\|\ol{\mb{v}}_{t-1}\|^2\right] +2n\beta^2\ol{\nu}^2\n\\
&+6L^2(1-\beta)^2\E\left[\|\x_t - \J\x_{t}\|^2 + \|\x_{t-1} - \J\x_{t-1}\|^2\right].   \label{vrb_sum}
\end{align}
\end{lem}

\begin{rmk}
Since~$\mb{v}_t$ is a conditionally biased estimator of~$\nabla \f(\x_t)$,~\eqref{vrb_a} and~\eqref{vrb_sum} do not directly imply each other and need to be established separately.
\end{rmk}
We emphasize that the contraction structure of the gradient variances shown in Lemma~\ref{vrb}
plays a crucial role in the convergence analysis.
The following contraction bounds on the consensus errors~$\E\left[\|\mb{x}_{t}-\J\mb{x}_{t}\|^2\right]$ are standard in decentralized algorithms based on gradient tracking, e.g.,~\cite{MP_Pu,GT-SARAH}; 
in particular, it follows directly from the~$\x_t$-update~\eqref{x_up}~and Young's inequality.  

\begin{lem}\label{consensus}
The following inequalities hold:~$\forall t\geq0$,
\begin{align}
\left\|\mb{x}_{t+1}-\J\mb{x}_{t+1}\right\|^2 \leq&~\dfrac{1+\lambda^2}{2}  \left\|\mb{x}_{t}-\J\mb{x}_{t}\right\|^2 \n\\
&+ \dfrac{2\alpha^2\lambda^2}{1-\lambda^2}\left\|\mb{y}_{t+1}-\J\mb{y}_{t+1}\right\|^2. \label{consensus_1} \\    
\left\|\mb{x}_{t+1}-\J\mb{x}_{t+1}\right\|^2 \leq&~2\lambda^2\left\|\mb{x}_{t}-\J\mb{x}_{t}\right\|^2  \n\\
&+ 2\alpha^2\lambda^2\left\|\mb{y}_{t+1}-\J\mb{y}_{t+1}\right\|^2.
\label{consensus_2}
\end{align}
\end{lem}
It is then clear from Lemma~\ref{consensus} that we need to further quantify the gradient tracking errors~$\E\left[\|\mb{y}_{t}-\J\mb{y}_{t}\|^2\right]$ in order to bound the consensus errors. These error bounds are shown in the following lemma. 

\begin{lem}\label{gt_bound}
We have the following.
\begin{enumerate}[(a)]
\item $\E\left[\|\mb{y}_1-\J\y_{1}\|^2\right]\leq \lambda^2\left\|\nabla\mb{f}\big(\x_{0}\big)\right\|^2 + \lambda^2n\ol{\nu}^2/b_0.$ \label{y1}
\item If~$0<\a\leq\frac{1-\lambda^2}{2\sqrt{42}\lambda^2L}$, then~$\forall t\geq1$,
\begin{align*}
&\E\left[\|\mb{y}_{t+1}-\J\mb{y}_{t+1}\|^2\right] \\
\leq&~\frac{3+\lambda^2}{4}\E\left[\|\mb{y}_{t}-\J\mb{y}_{t}\|^2\right] 
+ \frac{21\lambda^2nL^2\a^2}{1-\lambda^2}\E[\|\ol{\mb{v}}_{t-1}\|^2]
\n\\
&+ \frac{63\lambda^2L^2}{1-\lambda^2}\E\left[\|\x_{t-1}-\J\x_{t-1}\|^2\right] \n\\
&+\frac{7\lambda^2\beta^2}{1-\lambda^2}\E\left[\|\mb{v}_{t-1}
-\nabla\f(\x_{t-1})\|^2\right] 
+ 3\lambda^2n\beta^2\ol{\nu}^2. 
\end{align*} \label{yt}
\vspace*{-\baselineskip}
\end{enumerate}
\end{lem}

We note that establishing the contraction argument of gradient tracking errors in Lemma~\ref{gt_bound} requires a careful examination of the structure of the~$\mb{v}_t$-update. 

\subsection{Error Accumulations}

To proceed, we observe, from Lemma~\ref{vrb},~\ref{consensus}, and~\ref{gt_bound}, that the recursions of the gradient variances, consensus, and gradient tracking errors admit similar forms. Therefore, we abstract out formulas for the accumulation of the error recursions of this type in the following lemma. 

\begin{lem}\label{convs}
Let~$\{V_t\}_{t\geq0}$,~$\{R_t\}_{t\geq0}$ and~$\{Q_t\}_{t\geq0}$ be nonnegative sequences and~$C\geq0$ be some constant such that~$V_t\leq qV_{t-1} + qR_{t-1} + Q_t + C$,~$\forall t\geq1$, where~$q\in(0,1)$. Then the following inequality holds:~$\forall T\geq1$,
\begin{align}\label{convs1}
\sum_{t=0}^{T}V_t
\leq \frac{V_0}{1-q} 
+ \frac{1}{1-q}\sum_{t=0}^{T-1}R_{t} 
+ \frac{1}{1-q}\sum_{t=1}^{T}Q_{t} 
\!+ \frac{CT}{1-q}.
\end{align}
Similarly, if~$V_{t+1} \leq qV_{t} + R_{t-1} + C,\forall t\geq1$, then we have:~$\forall T\geq2$,
\begin{align}\label{convs2}
\sum_{t=1}^{T}V_{t} 
\leq \frac{V_{1}}{1-q} + \frac{1}{1-q}\sum_{t=0}^{T-2}R_t
+\frac{CT}{1-q}.
\end{align}
\end{lem}

Applying Lemma~\ref{convs} to Lemma~\ref{vrb} leads to the following upper bounds on the accumulated variances.
\begin{lem}\label{vrb_acc}
For any~$\beta\in(0,1)$, the following inequalities hold:~$\forall T\geq1$,
\begin{align}
&\sum_{t=0}^{T}\E\left[\|\ol{\mb{v}}_{t} - \ol{\nabla \f}(\x_t)\|^2\right] \n\\
\leq&~\frac{\ol{\nu}^2}{\beta b_0n}
+\frac{6L^2\a^2}{n\beta}\sum_{t=0}^{T-1}\E\left[\left\|\ol{\mb{v}}_{t}\right\|^2\right] 
\n\\
&+\frac{12L^2}{n^2\beta}
\sum_{t=0}^{T}\E\left[\|\x_{t} - \J\x_{t}\|^2\right]
+\frac{2\beta\ol{\nu}^2T}{n}, 
\label{vrb_a_acc}
\end{align}
and,~$\forall T\geq1$,
\begin{align}
&\sum_{t=0}^{T}\E\left[\|\mb{v}_t - \nabla\f(\x_t)\|^2\right] \n\\
\leq&~\frac{n\ol{\nu}^2}{\beta b_0}
+ \frac{6nL^2\a^2}{\beta}\sum_{t=0}^{T-1}\E\left[\|\ol{\mb{v}}_{t}\|^2\right] \n\\
&+ \frac{12L^2}{\beta}\sum_{t=0}^{T} \E\left[\|\x_{t} - \J\x_{t}\|^2\right]
+ 2n\beta\ol{\nu}^2T.
\label{vrb_s_acc}
\end{align}
\end{lem}

It can be observed that~\eqref{vrb_a_acc} in Lemma~\ref{vrb_acc} may be used to refine the left hand side of~\eqref{crt}. The remaining step, naturally, is to bound~$\sum_{t}\E\left[\|\mb{x}_{t}-\J\mb{x}_{t}\|^2\right]$ in terms of~$\sum_{t}\E\left[\|\ol{\mb{v}}_{t}\|^2\right]$. This result is provided in the following lemma that is obtained with the help of Lemma~\ref{consensus},~\ref{gt_bound},~\ref{convs}, and~\ref{vrb_acc}.

\begin{lem}\label{cons_bound_sum_final}
If~$0<\a\leq\frac{(1-\lambda^2)^2}{70\lambda^2L}$ and~$\beta\in(0,1)$, then the following inequality holds:~$\forall T\geq2$,
\begin{align*}
&\sum_{t=0}^{T}\frac{\E\left[\|\mb{x}_{t}-\J\mb{x}_{t}\|^2\right]}{n} 
\leq
\frac{2016\lambda^4L^2\a^4}{(1-\lambda^2)^4}\sum_{t=0}^{T-2}\E\left[\|\ol{\mb{v}}_{t}\|^2\right] \\
&+ \frac{32\lambda^4\a^2}{(1-\lambda^2)^3}\frac{\|\nabla\mb{f}\big(\x_{0}\big)\|^2}{n}
+ \left(\frac{7\beta}{1-\lambda^2}
+ 1\right)\frac{32\lambda^4\ol{\nu}^2\a^2}{(1-\lambda^2)^{3}b_0} \\
&+\left(\frac{14\beta}{1-\lambda^2}
+ 3\right)\frac{32\lambda^4\beta^2\ol{\nu}^2\a^2T}{(1-\lambda^2)^3}.
\end{align*}
\end{lem}
   
Finally, we note that Lemma~\ref{vrb_acc} and~\ref{cons_bound_sum_final} suffice to establish~\eqref{crt} and hence lead to Theorem~\ref{main0}; see the Appendix for details.

\begin{figure*}
\centering
\includegraphics[width=2.82in]{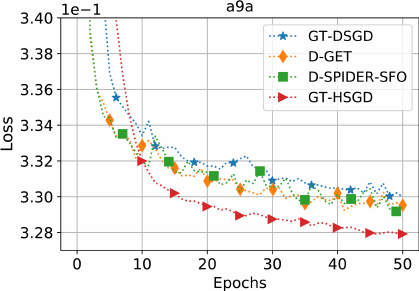}
\qquad
\includegraphics[width=2.82in]{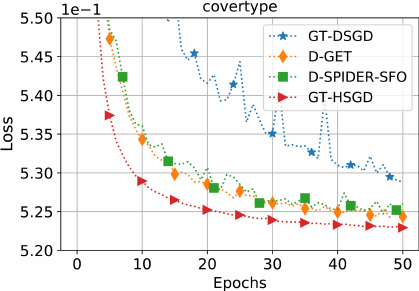}
\\
\vspace{0.2cm}
\includegraphics[width=2.82in]{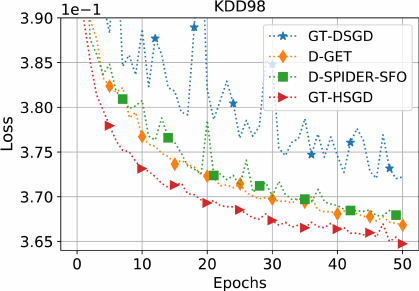}
\qquad\includegraphics[width=2.82in]{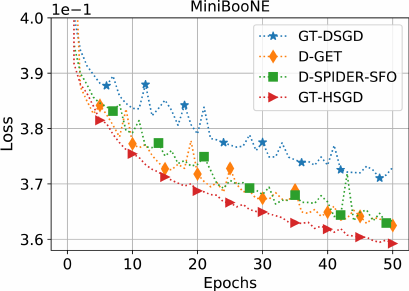}
\vspace*{-0.1cm} 
\caption{A comparison of \texttt{GT-HSGD} with other decentralized online stochastic gradient algorithms over the undirected exponential graph of~$20$ nodes on the a9a, covertype, KDD98, and MiniBooNE datasets.}
\label{alg_comp}
\end{figure*}

\section{Numerical Experiments}\label{sec_num_exp}
In this section, we illustrate our theoretical results on the convergence of the proposed \texttt{GT-HSGD} algorithm with the help of numerical experiments. 

\textbf{Model.} We consider a non-convex logistic regression model \cite{LR_NCVX} for decentralized binary classification. In particular, the decentralized non-convex optimization problem of interest takes the form~$\min_{\x\in\mathbb{R}^p}F(\mb{x})
:= \frac{1}{n}\sum_{i=1}^{n}f_{i}(\x)
+ r(\x),$
such that
\begin{align*}
f_{i}(\x) = \frac{1}{m}\sum_{j=1}^m\log\left[1+e^{-\langle\mb{x},\bds{\theta}_{ij}\rangle l_{ij}}\right]
\end{align*}
and $$r(\x) = R\sum_{k=1}^p\frac{[\x]_k^2}{1+[\x]_k^2},$$
where~$\bds{\theta}_{i,j}$ is the feature vector,~$l_{i,j}\in\{-1,+1\}$ is the corresponding binary label, and~$r(\x)$ is a non-convex regularizer. 
To simulate the online \SFO~setting described in Section~\ref{sec_problem_alg}, each node~$i$ is only able to \emph{sample with replacement} from its local data~$\{\bds{\theta}_{i,j},l_{i,j}\}_{j=1}^m$ and compute the corresponding (minibatch) stochastic gradient.  
Throughout all experiments, we set the number of the nodes to~$n = 20$ and the regularization parameter to~$R = 10^{-4}$. 

\textbf{Data.} To test the performance of the applicable decentralized algorithms, we distribute the a9a, covertype, KDD98, MiniBooNE datasets uniformly over the nodes and normalize the feature vectors such that~$\|\bds\theta_{i,j}\| = 1,\forall i,j$. The statistics of these datasets are provided in Table~\ref{datasets}.

\begin{table}[hbt]
\footnotesize
\renewcommand{\arraystretch}{1.5}
\vspace{-0.4cm}
\caption{Datasets used in numerical experiments, all available at \href{https://www.openml.org/}{https://www.openml.org/}.}
\begin{center}
\begin{tabular}{|c|c|c|c|}
\hline
\textbf{Dataset} & \textbf{train} ($nm$)  & \textbf{dimension} ($p$) \\ \hline
a9a & $48,\!840$ & $123$  \\ \hline
covertype & $100,\!000$ & $54$  \\ \hline
KDD98 &  $75,\!000$ & $477$ \\ \hline
MiniBooNE & $100,\!000$ & $11$  \\ \hline
\end{tabular}
\end{center}
\label{datasets}
\end{table}

\textbf{Network topology.} We consider the following network topologies: the undirected ring graph, the undirected and directed exponential graphs, and the complete graph; see \cite{tutorial_nedich,PIEEE_Xin,SGP_ICML,DSGD_NIPS} for detailed configurations of these graphs. For all graphs, the associated doubly stochastic weights are set to be equal. The resulting second largest singular value~$\lambda$ of the weight matrices are~$0.98,0.75,0.67,0$, respectively, demonstrating a significant difference in the algebraic connectivity of these graphs. 

\textbf{Performance measure.}
We measure the performance of the decentralized algorithms in question by the decrease of the global cost function value~$F(\ol{\x})$, to which we refer as loss, versus epochs, where~$\ol{\x} = \frac{1}{n}\sum_{i=1}^n\x_i$ with $\x_i$ being the model at node~$i$ and each epoch contains~$m$ stochastic gradient computations at each node. 

\subsection{Comparison with the Existing Decentralized Stochastic Gradient Methods.} \label{sec_num_comp}
We conduct a performance comparison of \texttt{GT-HSGD} with \texttt{GT-DSGD} \cite{MP_Pu,GNSD,improved_DSGT_Xin}, \texttt{D-GET} \cite{D_Get}, and \texttt{D-SPIDER-SFO} \cite{D-SPIDER-SFO} over the undirected exponential graph of~$20$ nodes. Note that we use \texttt{GT-DSGD} to represent methods that do not incorporate online variance reduction techniques, since it
in general matches or outperforms \texttt{DSGD} \cite{DSGD_NIPS} and has a similar performance with \texttt{D2} \cite{D2} and \texttt{D-PD-SGD} \cite{PD_SGD}. 

\textbf{Parameter tuning.} We set the parameters of \texttt{GT-HSGD}, \texttt{GT-DSGD}, \texttt{D-GET}, and \texttt{D-SPIDER-SFO} according to the following procedures. \emph{First}, we find a very large step-size candidate set for each algorithm in comparison. \emph{Second}, we choose the minibatch size candidate set for all algorithms as $\mc{B} := \{1,4,8,16,32,64,128,256,512,1024\}$: the minibatch size of \texttt{GT-DSGD}, the minibatch size of \texttt{GT-HSGD} at $t=0$, the minibatch size of \texttt{D-GET} and \texttt{D-SPIDER-SFO} at inner- and outer-loop are all chosen from~$\mc{B}$. \emph{Third}, for \texttt{D-GET} and \texttt{D-SPIDER-SFO}, we choose the inner-loop length candidate set as $\{\frac{m}{20b},\frac{m}{19b},\cdots,\frac{m}{b},\frac{2m}{b},\cdots,\frac{20m}{b}\}$, where~$m$ is the local data size and~$b$ is the minibatch size at the inner-loop. \emph{Fourth}, we iterate over all combinations of parameters for each algorithm to find its best performance. In particular, we find that the best performance of \texttt{GT-HSGD} is attained with a small~$\beta$ and a relatively large~$\alpha$ as Corollary~\ref{main} suggests.

The experimental results are provided in Fig.~\ref{alg_comp}, where we observe that \texttt{GT-HSGD} achieves faster convergence than the other algorithms in comparison on those four datasets. This observation is coherent with our main convergence results that \texttt{GT-HSGD} achieves a lower oracle complexity than the existing approaches; see Table~\ref{comp_oc}. 


\begin{figure}
\centering
\includegraphics[width=2.82in]{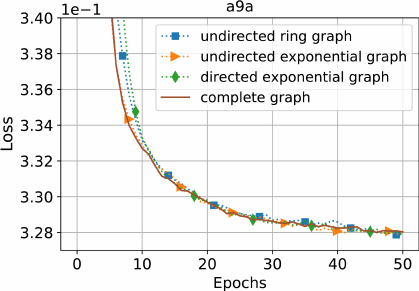}\\
\vspace{0.2cm}
\includegraphics[width=2.82in]{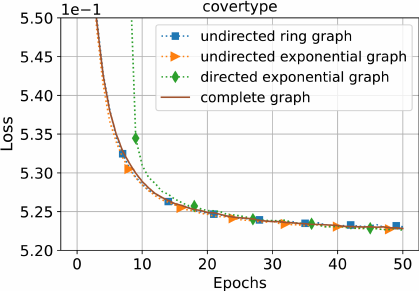}
\caption{Convergence behaviors of \texttt{GT-HSGD} over different network topologies on the a9a and covertype datasets.}
\label{net_indep}
\vspace{-0.4cm}
\end{figure}

\subsection{Topology-Independent Rate of \texttt{GT-HSGD}}
We test the performance of \texttt{GT-HSGD} over different network topologies. 
In particular, we follow the procedures described in Section~\ref{sec_num_comp} to find the best set of parameters for \texttt{GT-HSGD} over the \emph{complete graph} and then use this parameter set for other graphs. 
The corresponding experimental results are presented in Fig.~\ref{net_indep}. Clearly, it can be observed that when the number of iterations is large enough, that is to say, the required error tolerance is small enough, the convergence rate of \texttt{GT-HSGD} is not affected by the underlying network topology. This interesting phenomenon is consistent with our convergence theory; see Corollary~\ref{main} and the related discussion in Section~\ref{sec_main_results}. 

\vspace{-0.1cm}
\section{Conclusion}\label{sec_conc}
In this paper, we investigate decentralized stochastic optimization to minimize a sum of smooth non-convex cost functions over a network of nodes. Assuming that each node has access to a stochastic first-order oracle, we propose \texttt{GT-HSGD}, a novel single-loop decentralized algorithm that leverages local hybrid variance reduction and gradient tracking to achieve provably fast convergence and robust performance. Compared with the existing online variance-reduced methods, \texttt{GT-HSGD} achieves a lower oracle complexity with a more practical implementation. We further show that \texttt{GT-HSGD} achieves a network topology-independent oracle complexity, when the required error tolerance is small enough, leading to a linear speedup with respect to the centralized optimal methods that execute on a single node.

\vspace{-0.1cm}
\section*{Acknowledgments}
The work of Ran Xin and Soummya Kar was supported in part by NSF under Award \#1513936. The work
of Usman A. Khan was supported in part by NSF under Awards \#1903972 and
\#1935555.

\bibliography{ICML_xin}

\begin{thebibliography}{51}
\providecommand{\natexlab}[1]{#1}
\providecommand{\url}[1]{\texttt{#1}}
\expandafter\ifx\csname urlstyle\endcsname\relax
  \providecommand{\doi}[1]{doi: #1}\else
  \providecommand{\doi}{doi: \begingroup \urlstyle{rm}\Url}\fi

\bibitem[Alghunaim et~al.(2020)Alghunaim, Ryu, Yuan, and Sayed]{DGT_NIPS}
Alghunaim, S.~A., Ryu, E., Yuan, K., and Sayed, A.~H.
\newblock Decentralized proximal gradient algorithms with linear convergence
  rates.
\newblock \emph{IEEE Trans. Autom. Control}, 2020.

\bibitem[Antoniadis et~al.(2011)Antoniadis, Gijbels, and Nikolova]{LR_NCVX}
Antoniadis, A., Gijbels, I., and Nikolova, M.
\newblock Penalized likelihood regression for generalized linear models with
  non-quadratic penalties.
\newblock \emph{Annals of the Institute of Statistical Mathematics},
  63\penalty0 (3):\penalty0 585--615, 2011.

\bibitem[Arjevani et~al.(2019)Arjevani, Carmon, Duchi, Foster, Srebro, and
  Woodworth]{lowerbound_sgd}
Arjevani, Y., Carmon, Y., Duchi, J.~C., Foster, D.~J., Srebro, N., and
  Woodworth, B.
\newblock Lower bounds for non-convex stochastic optimization.
\newblock \emph{arXiv preprint arXiv:1912.02365}, 2019.

\bibitem[Assran et~al.(2019)Assran, Loizou, Ballas, and Rabbat]{SGP_ICML}
Assran, M., Loizou, N., Ballas, N., and Rabbat, M.
\newblock Stochastic gradient push for distributed deep learning.
\newblock In \emph{Proceedings of the 36th International Conference on Machine
  Learning}, pp.\  97: 344--353, 2019.

\bibitem[Chen \& Sayed(2015)Chen and Sayed]{DSGD1_Chen}
Chen, J. and Sayed, A.~H.
\newblock On the learning behavior of adaptive networks—part i: Transient
  analysis.
\newblock \emph{IEEE Transactions on Information Theory}, 61\penalty0
  (6):\penalty0 3487--3517, 2015.

\bibitem[Cutkosky \& Orabona(2019)Cutkosky and Orabona]{mSARAH}
Cutkosky, A. and Orabona, F.
\newblock Momentum-based variance reduction in non-convex sgd.
\newblock In \emph{Adv. Neural Inf. Process. Syst.}, pp.\  15236--15245, 2019.

\bibitem[Di~Lorenzo \& Scutari(2016)Di~Lorenzo and Scutari]{NEXT_scutari}
Di~Lorenzo, P. and Scutari, G.
\newblock {NEXT: I}n-network nonconvex optimization.
\newblock \emph{IEEE Trans. Signal Inf. Process. Netw. Process.}, 2\penalty0
  (2):\penalty0 120--136, 2016.

\bibitem[Fang et~al.(2018)Fang, Li, Lin, and Zhang]{spider}
Fang, C., Li, C.~J., Lin, Z., and Zhang, T.
\newblock {SPIDER:} near-optimal non-convex optimization via stochastic
  path-integrated differential estimator.
\newblock In \emph{Proc. Adv. Neural Inf. Process. Syst.}, pp.\  689--699,
  2018.

\bibitem[Jakoveti{\'c}(2018)]{GT_jakovetic}
Jakoveti{\'c}, D.
\newblock A unification and generalization of exact distributed first-order
  methods.
\newblock \emph{IEEE Trans. Signal Inf. Process. Netw. Process.}, 5\penalty0
  (1):\penalty0 31--46, 2018.

\bibitem[Kar et~al.(2012)Kar, Moura, and Ramanan]{DGD_Kar}
Kar, S., Moura, J. M.~F., and Ramanan, K.
\newblock Distributed parameter estimation in sensor networks: Nonlinear
  observation models and imperfect communication.
\newblock \emph{IEEE Trans. Inf. Theory}, 58\penalty0 (6):\penalty0 3575--3605,
  2012.

\bibitem[Li et~al.(2020{\natexlab{a}})Li, Cen, Chen, and Chi]{Network-DANE}
Li, B., Cen, S., Chen, Y., and Chi, Y.
\newblock Communication-efficient distributed optimization in networks with
  gradient tracking and variance reduction.
\newblock \emph{J. Mach. Learn. Res.}, 21\penalty0 (180):\penalty0 1--51,
  2020{\natexlab{a}}.

\bibitem[Li \& Lin(2020)Li and Lin]{EXTRA_revisit}
Li, H. and Lin, Z.
\newblock Revisiting extra for smooth distributed optimization.
\newblock \emph{SIAM J. Optim.}, 30\penalty0 (3):\penalty0 1795--1821, 2020.

\bibitem[Li et~al.(2020{\natexlab{b}})Li, Lin, and Fang]{AccDVR}
Li, H., Lin, Z., and Fang, Y.
\newblock Optimal accelerated variance reduced {EXTRA} and {DIGing} for
  strongly convex and smooth decentralized optimization.
\newblock \emph{arXiv preprint arXiv:2009.04373}, 2020{\natexlab{b}}.

\bibitem[Li et~al.(2019)Li, Shi, and Yan]{NIDS}
Li, Z., Shi, W., and Yan, M.
\newblock A decentralized proximal-gradient method with network independent
  step-sizes and separated convergence rates.
\newblock \emph{IEEE Trans. Signal Process.}, 67\penalty0 (17):\penalty0
  4494--4506, 2019.

\bibitem[Lian et~al.(2017)Lian, Zhang, Zhang, Hsieh, Zhang, and Liu]{DSGD_NIPS}
Lian, X., Zhang, C., Zhang, H., Hsieh, C.-J., Zhang, W., and Liu, J.
\newblock {Can decentralized algorithms outperform centralized algorithms? {A}
  case study for decentralized parallel stochastic gradient descent}.
\newblock In \emph{Adv. Neural Inf. Process. Syst.}, pp.\  5330--5340, 2017.

\bibitem[Liu et~al.(2020)Liu, Nguyen, and Tran-Dinh]{HSARAH}
Liu, D., Nguyen, L.~M., and Tran-Dinh, Q.
\newblock An optimal hybrid variance-reduced algorithm for stochastic composite
  nonconvex optimization.
\newblock \emph{arXiv preprint arXiv:2008.09055}, 2020.

\bibitem[L{\"u} et~al.(2020)L{\"u}, Liao, Li, and Huang]{DVR_LHQ}
L{\"u}, Q., Liao, X., Li, H., and Huang, T.
\newblock A computation-efficient decentralized algorithm for composite
  constrained optimization.
\newblock \emph{IEEE Trans. Signal Inf. Process. Netw.}, 6:\penalty0 774--789,
  2020.

\bibitem[Lu et~al.(2019)Lu, Zhang, Sun, and Hong]{GNSD}
Lu, S., Zhang, X., Sun, H., and Hong, M.
\newblock {GNSD: A gradient-tracking based nonconvex stochastic algorithm for
  decentralized optimization}.
\newblock In \emph{2019 IEEE Data Science Workshop}, pp.\  315--321, 2019.

\bibitem[Mokhtari \& Ribeiro(2016)Mokhtari and Ribeiro]{DSA}
Mokhtari, A. and Ribeiro, A.
\newblock {DSA: Decentralized double stochastic averaging gradient algorithm}.
\newblock \emph{J. Mach. Learn. Res.}, 17\penalty0 (1):\penalty0 2165--2199,
  2016.

\bibitem[Nedi\'{c} \& Ozdaglar(2009)Nedi\'{c} and Ozdaglar]{DGD_nedich}
Nedi\'{c}, A. and Ozdaglar, A.
\newblock Distributed subgradient methods for multi-agent optimization.
\newblock \emph{IEEE Trans. Autom. Control}, 54\penalty0 (1):\penalty0 48,
  2009.

\bibitem[Nedi{\'c} et~al.(2018)Nedi{\'c}, Olshevsky, and
  Rabbat]{tutorial_nedich}
Nedi{\'c}, A., Olshevsky, A., and Rabbat, M.~G.
\newblock Network topology and communication-computation tradeoffs in
  decentralized optimization.
\newblock \emph{P. IEEE}, 106\penalty0 (5):\penalty0 953--976, 2018.

\bibitem[Nedich et~al.(2017)Nedich, Olshevsky, and Shi]{DIGing}
Nedich, A., Olshevsky, A., and Shi, W.
\newblock Achieving geometric convergence for distributed optimization over
  time-varying graphs.
\newblock \emph{SIAM J. Optim.}, 27\penalty0 (4):\penalty0 2597--2633, 2017.

\bibitem[Nemirovski et~al.(2009)Nemirovski, Juditsky, Lan, and
  Shapiro]{SGD_Nemirovski}
Nemirovski, A., Juditsky, A., Lan, G., and Shapiro, A.
\newblock Robust stochastic approximation approach to stochastic programming.
\newblock \emph{SIAM J. Optim.}, 19\penalty0 (4):\penalty0 1574--1609, 2009.

\bibitem[Nesterov(2018)]{nesterov_book}
Nesterov, Y.
\newblock \emph{Lectures on convex optimization}, volume 137.
\newblock Springer, 2018.

\bibitem[Nguyen et~al.(2017)Nguyen, Liu, Scheinberg, and Takac]{SARAH}
Nguyen, L.~M., Liu, J., Scheinberg, K., and Takac, M.
\newblock {SARAH: A} novel method for machine learning problems using
  stochastic recursive gradient.
\newblock In \emph{Proc. 34th Int. Conf. Mach. Learn.}, pp.\  2613--2621, 2017.

\bibitem[Pan et~al.(2020)Pan, Liu, and Wang]{D-SPIDER-SFO}
Pan, T., Liu, J., and Wang, J.
\newblock {D-SPIDER-SFO: A} decentralized optimization algorithm with faster
  convergence rate for nonconvex problems.
\newblock In \emph{Proceedings of the AAAI Conference on Artificial
  Intelligence}, volume~34, pp.\  1619--1626, 2020.

\bibitem[Pham et~al.(2020)Pham, Nguyen, Phan, and Tran-Dinh]{sarah_ncvx}
Pham, N.~H., Nguyen, L.~M., Phan, D.~T., and Tran-Dinh, Q.
\newblock {ProxSARAH:} an efficient algorithmic framework for stochastic
  composite nonconvex optimization.
\newblock \emph{J. Mach. Learn. Res.}, 21\penalty0 (110):\penalty0 1--48, 2020.

\bibitem[Pu \& Nedich(2020)Pu and Nedich]{MP_Pu}
Pu, S. and Nedich, A.
\newblock Distributed stochastic gradient tracking methods.
\newblock \emph{Math. Program.}, pp.\  1--49, 2020.

\bibitem[Qu \& Li(2017)Qu and Li]{harnessing}
Qu, G. and Li, N.
\newblock Harnessing smoothness to accelerate distributed optimization.
\newblock \emph{IEEE Trans. Control. Netw. Syst.}, 5\penalty0 (3):\penalty0
  1245--1260, 2017.

\bibitem[Rajawat \& Kumar(2020)Rajawat and Kumar]{DPD_VR}
Rajawat, K. and Kumar, C.
\newblock A primal-dual framework for decentralized stochastic optimization.
\newblock \emph{arXiv preprint arXiv:2012.04402}, 2020.

\bibitem[Shi et~al.(2015)Shi, Ling, Wu, and Yin]{EXTRA}
Shi, W., Ling, Q., Wu, G., and Yin, W.
\newblock {EXTRA: A}n exact first-order algorithm for decentralized consensus
  optimization.
\newblock \emph{SIAM J. Optim.}, 25\penalty0 (2):\penalty0 944--966, 2015.

\bibitem[Sun et~al.(2020)Sun, Lu, and Hong]{D_Get}
Sun, H., Lu, S., and Hong, M.
\newblock Improving the sample and communication complexity for decentralized
  non-convex optimization: Joint gradient estimation and tracking.
\newblock In \emph{Proceedings of the 37th International Conference on Machine
  Learning}, volume 119, pp.\  9217--9228, 13--18 Jul 2020.

\bibitem[Taheri et~al.(2020)Taheri, Mokhtari, Hassani, and Pedarsani]{QSGP}
Taheri, H., Mokhtari, A., Hassani, H., and Pedarsani, R.
\newblock Quantized decentralized stochastic learning over directed graphs.
\newblock In \emph{International Conference on Machine Learning}, pp.\
  9324--9333, 2020.

\bibitem[Tang et~al.(2018)Tang, Lian, Yan, Zhang, and Liu]{D2}
Tang, H., Lian, X., Yan, M., Zhang, C., and Liu, J.
\newblock {$D^{2}$: D}ecentralized training over decentralized data.
\newblock In \emph{International Conference on Machine Learning}, pp.\
  4848--4856, 2018.

\bibitem[Tran-Dinh et~al.(2020)Tran-Dinh, Pham, Phan, and Nguyen]{HSARAH_0}
Tran-Dinh, Q., Pham, N.~H., Phan, D.~T., and Nguyen, L.~M.
\newblock A hybrid stochastic optimization framework for stochastic composite
  nonconvex optimization.
\newblock \emph{Math. Program.}, 2020.

\bibitem[Tsitsiklis et~al.(1986)Tsitsiklis, Bertsekas, and
  Athans]{DGD_tsitsiklis}
Tsitsiklis, J., Bertsekas, D., and Athans, M.
\newblock Distributed asynchronous deterministic and stochastic gradient
  optimization algorithms.
\newblock \emph{IEEE Trans. Autom. Control}, 31\penalty0 (9):\penalty0
  803--812, 1986.

\bibitem[Vlaski \& Sayed(2019)Vlaski and Sayed]{DSGD_vlaski_2}
Vlaski, S. and Sayed, A.~H.
\newblock Distributed learning in non-convex environments--{Part II:
  P}olynomial escape from saddle-points.
\newblock \emph{arXiv:1907.01849}, 2019.

\bibitem[Wang et~al.(2019)Wang, Ji, Zhou, Liang, and Tarokh]{spiderboost}
Wang, Z., Ji, K., Zhou, Y., Liang, Y., and Tarokh, V.
\newblock Spiderboost and momentum: Faster variance reduction algorithms.
\newblock In \emph{Proc. Adv. Neural Inf. Process. Syst.}, pp.\  2403--2413,
  2019.

\bibitem[Xi et~al.(2017)Xi, Xin, and Khan]{add-opt}
Xi, C., Xin, R., and Khan, U.~A.
\newblock {ADD-OPT: Accelerated distributed directed optimization}.
\newblock \emph{IEEE Trans. Autom. Control}, 63\penalty0 (5):\penalty0
  1329--1339, 2017.

\bibitem[Xin et~al.(2020{\natexlab{a}})Xin, Kar, and Khan]{SPM_Xin}
Xin, R., Kar, S., and Khan, U.~A.
\newblock Decentralized stochastic optimization and machine learning: A unified
  variance-reduction framework for robust performance and fast convergence.
\newblock \emph{IEEE Signal Process. Mag.}, 37\penalty0 (3):\penalty0 102--113,
  2020{\natexlab{a}}.

\bibitem[Xin et~al.(2020{\natexlab{b}})Xin, Khan, and Kar]{GT-SARAH}
Xin, R., Khan, U.~A., and Kar, S.
\newblock A near-optimal stochastic gradient method for decentralized
  non-convex finite-sum optimization.
\newblock \emph{arXiv preprint arXiv:2008.07428}, 2020{\natexlab{b}}.

\bibitem[Xin et~al.(2020{\natexlab{c}})Xin, Khan, and Kar]{GTSAGA_NCVX}
Xin, R., Khan, U.~A., and Kar, S.
\newblock A fast randomized incremental gradient method for decentralized
  non-convex optimization.
\newblock \emph{arXiv preprint arXiv:2011.03853}, 2020{\natexlab{c}}.

\bibitem[Xin et~al.(2020{\natexlab{d}})Xin, Khan, and Kar]{GTVR}
Xin, R., Khan, U.~A., and Kar, S.
\newblock Variance-reduced decentralized stochastic optimization with
  accelerated convergence.
\newblock \emph{IEEE Trans. Signal Process.}, 68:\penalty0 6255--6271,
  2020{\natexlab{d}}.

\bibitem[Xin et~al.(2020{\natexlab{e}})Xin, Khan, and Kar]{improved_DSGT_Xin}
Xin, R., Khan, U.~A., and Kar, S.
\newblock An improved convergence analysis for decentralized online stochastic
  non-convex optimization.
\newblock \emph{arXiv preprint arXiv:2008.04195}, 2020{\natexlab{e}}.

\bibitem[Xin et~al.(2020{\natexlab{f}})Xin, Pu, Nedi{\'c}, and Khan]{PIEEE_Xin}
Xin, R., Pu, S., Nedi{\'c}, A., and Khan, U.~A.
\newblock A general framework for decentralized optimization with first-order
  methods.
\newblock \emph{Proceedings of the IEEE}, 108\penalty0 (11):\penalty0
  1869--1889, 2020{\natexlab{f}}.

\bibitem[Xu et~al.(2015)Xu, Zhu, Soh, and Xie]{GT_CDC}
Xu, J., Zhu, S., Soh, Y.~C., and Xie, L.
\newblock Augmented distributed gradient methods for multi-agent optimization
  under uncoordinated constant stepsizes.
\newblock In \emph{Proc. IEEE Conf. Decis. Control}, pp.\  2055--2060, 2015.

\bibitem[Xu et~al.(2020)Xu, Tian, Sun, and Scutari]{PD_Xu}
Xu, J., Tian, Y., Sun, Y., and Scutari, G.
\newblock Distributed algorithms for composite optimization: Unified and tight
  convergence analysis.
\newblock \emph{arXiv:2002.11534}, 2020.

\bibitem[Yi et~al.(2020)Yi, Zhang, Yang, Chai, and Johansson]{PD_SGD}
Yi, X., Zhang, S., Yang, T., Chai, T., and Johansson, K.~H.
\newblock A primal-dual {SGD} algorithm for distributed nonconvex optimization.
\newblock \emph{arXiv preprint arXiv:2006.03474}, 2020.

\bibitem[Yuan et~al.(2018)Yuan, Ying, Liu, and Sayed]{DAVRG}
Yuan, K., Ying, B., Liu, J., and Sayed, A.~H.
\newblock Variance-reduced stochastic learning by~networked agents under random
  reshuffling.
\newblock \emph{IEEE Trans. Signal Process.}, \!67\penalty0 (2):\penalty0
  351--366, 2018.

\bibitem[Yuan et~al.(2020)Yuan, Alghunaim, Ying, and Sayed]{SED}
Yuan, K., Alghunaim, S.~A., Ying, B., and Sayed, A.~H.
\newblock On the influence of bias-correction on distributed stochastic
  optimization.
\newblock \emph{IEEE Trans. Signal Process.}, 2020.

\bibitem[Zhou et~al.(2020)Zhou, Xu, and Gu]{SNVRG}
Zhou, D., Xu, P., and Gu, Q.
\newblock Stochastic nested variance reduction for nonconvex optimization.
\newblock \emph{J. Mach. Learn. Res.}, 2020.

\end{thebibliography}
\bibliographystyle{icml2021}
\newpage

\onecolumn
\appendix
\section{Proof of Lemma~\ref{DS0}}\label{proof_DS0}
We recall the standard Descent Lemma~\cite{nesterov_book}, i.e.,~$\forall\mb{x},\mb{y}\in\mathbb{R}^p$,
\begin{equation}\label{DL}
F(\mb y) \leq F(\mb x) + \langle \nabla F(\mb x), \mb{y} - \mb{x}\rangle
+ \frac{L}{2}\left\|\mb y - \mb x\right\|^2,
\end{equation}
since the global function~$F$ is~$L$-smooth.
Setting~$\y=\ol{\x}_{t+1}$ and~$\x=\ol{\x}_{t}$ in~\eqref{DL} and using~\eqref{ave} , we have:~$\forall t\geq0$,
\begin{align}\label{DS000}
F(\ol{\x}_{t+1}) 
\leq&~F(\ol{\x}_{t}) - \big\langle\nabla F(\ol{\x}_{t}),\ol{\x}_{t+1} - \ol{\x}_{t}\big\rangle
+ \frac{L}{2}\left\|\ol{\x}_{t+1} - \ol{\x}_{t}\right\|^2 \n\\
\leq&~F(\ol{\x}_{t}) - \alpha\big\langle\nabla F(\ol{\x}_{t}),\ol{\mb{v}}_{t}\big\rangle
+ \frac{L\alpha^2 }{2}\left\|\ol{\mb{v}}_{t}\right\|^2.
\end{align}
Using~$\langle \mb{a},\mb{b} \rangle = \frac{1}{2}\left(\|\mb a\|^2 + \|\mb b\|^2 - \|\mb{a}-\mb{b}\|^2\right), \forall\mb{a},\mb{b}\in\mbb{R}^p$, in~\eqref{DS000} gives: for~$0<\a\leq\frac{1}{2L}$ and~$\forall t\geq0$,
\begin{align}\label{DS00}
F(\ol{\x}_{t+1}) \leq&~F(\ol{\x}_{t}) 
- \frac{\a}{2}\left\|\nabla F(\ol{\x}_{t})\right\|^2
- \left(\frac{\a}{2}-\frac{L\a^2}{2}\right)\left\|\ol{\mb{v}}_{t}\right\|^2 + \frac{\a}{2}\left\|\ol{\mb{v}}_{t}-\nabla F(\ol{\x}_{t})\right\|^2, \nonumber\\
\leq&~F(\ol{\x}_{t}) - \frac{\a}{2}\left\|\nabla F(\ol{\x}_{t})\right\|^2
- \left(\frac{\a}{2}-\frac{L\a^2}{2}\right)\left\|\ol{\mb{v}}_{t}\right\|^2 + \a\left\|\ol{\mb{v}}_{t}-\ol{\nabla\mb{f}}(\x_{t})\right\|^2 + \a\left\|\ol{\nabla\mb{f}}(\x_{t}) - \nabla F(\ol{\x}_{t})\right\|^2, \nonumber\\
\stackrel{(i)}{\leq}&~F(\ol{\x}_{t}) - \frac{\a}{2}\left\|\nabla F(\ol{\x}_{t})\right\|^2
- \frac{\a}{4}\left\|\ol{\mb{v}}_{t}\right\|^2 + \a\left\|\ol{\mb{v}}_{t}-\ol{\nabla\mb{f}}(\x_{t})\right\|^2 
+ \frac{\a L^2}{n}\left\|\mb{x}_{t}-\J\x_{t}\right\|^2,
\end{align}
where~$(i)$ is due to Lemma~\ref{basic}(\ref{Lbound}) and that~$\frac{L\a^2}{2}\leq\frac{\a}{4}$ since~$0<\a\leq\frac{1}{2L}$. Rearranging~\eqref{DS00}, we have: for~$0<\a\leq\frac{1}{2L}$ and~$\forall t\geq0$,
\begin{align}\label{DS01}
\left\|\nabla F(\ol{\x}_{t})\right\|^2 
\leq 
\frac{2(F(\ol{\x}_{t}) - F(\ol{\x}_{t+1}))}{\a}
- \frac{1}{2}\left\|\ol{\mb{v}}_{t}\right\|^2 + 2\left\|\ol{\mb{v}}_{t}-\ol{\nabla\mb{f}}(\x_{t})\right\|^2 
+ \frac{2L^2}{n}\left\|\mb{x}_{t}-\J\x_{t}\right\|^2.
\end{align}
Taking the telescoping sum of~\eqref{DS01} over~$t$ from~$0$ to~$T$,~$\forall T\geq0$ and using the fact that~$F$ bounded below by~$F^*$ in the resulting inequality finishes the proof.

\section{Proof of Lemma~\ref{vrb}}\label{proof_vrb}
\subsection{Proof of Eq.~\eqref{vrb_a}}
We recall that the update of each local stochastic gradient estimator~$\mb{v}_t^i,\forall t\geq1$, in~\eqref{vi} may be written equivalently as follows:
\begin{align*}
\mb{v}_{t}^i 
= \beta\g_i(\x_t^i,\X_{t}^i) 
+ (1-\beta)\Big(\g_i(\x_t^i,\X_{t}^i) - \g_i(\x_{t-1}^i,\X_{t}^i) + \mb{v}_{t-1}^i\Big),
\end{align*}
where~$\beta\in(0,1)$. We have:~$\forall t\geq1$ and~$\forall i\in\mc{V}$,
\begin{align}\label{vr_b0}
\mb{v}_{t}^i - \nabla f_i(\x_t^i)
=&~\beta\g_i(\x_t^i,\X_{t}^i) 
+ (1-\beta)\Big(\g_i(\x_t^i,\X_{t}^i) - \g_i(\x_{t-1}^i,\X_{t}^i) + \mb{v}_{t-1}^i\Big)
- \beta \nabla f_i(\x_t^i) - (1 - \beta) \nabla f_i(\x_t^i) \n\\
=&~\beta\Big(\g_i(\x_t^i,\X_{t}^i) - \nabla f_i(\x_t^i)\Big)
+ (1-\beta)\Big(\g_i(\x_t^i,\X_{t}^i) - \g_i(\x_{t-1}^i,\X_{t}^i) + \mb{v}_{t-1}^i
- \nabla f_i(\x_t^i)\Big) \n\\
=&~\beta\Big(\g_i(\x_t^i,\X_{t}^i) - \nabla f_i(\x_t^i)\Big)
+ (1-\beta)\Big(\g_i(\x_t^i,\X_{t}^i) - \g_i(\x_{t-1}^i,\X_{t}^i) + \nabla f_i(\x_{t-1}^i)
- \nabla f_i(\x_t^i)\Big) \n\\
&+(1-\beta)\Big(\mb{v}_{t-1}^i - \nabla f_i(\x_{t-1}^i)\Big).
\end{align}
In~\eqref{vr_b0}, we observe that~$\forall t\geq1$ and~$\forall i\in\mc{V}$,
\begin{align}
&\E\Big[\g_i(\x_t^i,\X_{t}^i) - \nabla f_i(\x_t^i)|\F_t\Big] = \mb{0}_p, \label{O1}
\\
&\E\Big[\g_i(\x_t^i,\X_{t}^i) - \g_i(\x_{t-1}^i,\X_{t}^i) + \nabla f_i(\x_{t-1}^i)
- \nabla f_i(\x_t^i)|\F_t\Big] = \mb{0}_p,    \label{O2}
\end{align}
by the definition of the filtration~$\F_t$ in~\eqref{Ft}.
Averaging~\eqref{vr_b0} over~$i$ from~$1$ to~$n$ gives:~$\forall t\geq0$,
\begin{align}\label{vr_b_a1}
\ol{\mb{v}}_{t} - \ol{\nabla \f}(\x_t)
=&~(1-\beta)\Big(\ol{\mb{v}}_{t-1} - \ol{\nabla \f}(\x_{t-1})\Big)
+\beta\cdot\underbrace{\frac{1}{n}\sum_{i=1}^n\Big(\g_i(\x_t^i,\X_{t}^i) - \nabla f_i(\x_t^i)\Big)}_{=:\mb{s}_t} \n\\
&+ (1-\beta)\cdot\underbrace{\frac{1}{n}\sum_{i=1}^n\Big(\g_i(\x_t^i,\X_{t}^i) - \g_i(\x_{t-1}^i,\X_{t}^i) + \nabla f_i(\x_{t-1}^i)
-\nabla f_i(\x_t^i)\Big)}_{=:\mb{z}_t}.
\end{align}
Note that~$\E[\mb{s}_t|\F_t] = \E[\mb{z}_t|\F_t] = \mb{0}_p$ by~\eqref{O1} and~\eqref{O2}. In light of~\eqref{vr_b_a1}, we have:~$\forall t\geq1$,
\begin{align}\label{vr_b_a2}
\E\Big[\|\ol{\mb{v}}_{t} - \ol{\nabla \f}(\x_t)\|^2|\F_t\Big]  =&~(1-\beta)^2\|\ol{\mb{v}}_{t-1} - \ol{\nabla \f}(\x_{t-1})\|^2
+\E\left[\left\|\beta\s_t + (1-\beta)\z_t\right\|^2|\F_t\right] \n\\
&+2\E\Big[\Big\langle(1-\beta)\left(\ol{\mb{v}}_{t-1} - \ol{\nabla \f}(\x_{t-1})\right),\beta\s_t + (1-\beta)\z_t\Big\rangle|\F_t\Big]\n\\
\stackrel{(i)}{=}&~(1-\beta)^2\|\ol{\mb{v}}_{t-1} - \ol{\nabla \f}(\x_{t-1})\|^2
+\E\left[\left\|\beta\s_t + (1-\beta)\z_t\right\|^2|\F_t\right] \n\\
\leq&~(1-\beta)^2\|\ol{\mb{v}}_{t-1} - \ol{\nabla \f}(\x_{t-1})\|^2
+2\beta^2\E\left[\left\|\s_t\right\|^2|\F_t\right] 
+2(1-\beta)^2\E\left[\left\|\z_t\right\|^2|\F_t\right],
\end{align}
where~$(i)$ is due to
$$\E\Big[\Big\langle(1-\beta)\left(\ol{\mb{v}}_{t-1} - \ol{\nabla \f}(\x_{t-1})\right),\beta\s_t + (1-\beta)\z_t\Big\rangle |\F_t\Big]=0,$$
since~$\E[\s_t|\F_t] = \E[\z_t|\F_t] = \mb{0}_p$ and~$(\ol{\mb{v}}_{t-1} - \ol{\nabla \f}(\x_{t-1}))$ is~$\F_t$-measurable. We next bound the second and the third term in~\eqref{vr_b_a2} respectively. 
For the second term in~\eqref{vr_b_a2}, we observe that~$\forall t\geq1$,
\begin{align}\label{st}
\E\left[\|\s_t\|^2\right]        
=&~\frac{1}{n^2}\sum_{i=1}^n\E\left[\left\|\g_i(\x_t^i,\X_{t}^i) - \nabla f_i(\x_t^i)\right\|^2\right] 
+\frac{1}{n^2}\sum_{i\neq j}\E\left[\Big\langle\g_i(\x_t^i,\X_{t}^i) - \nabla f_i(\x_t^i),\g_j(\x_t^j,\X_{t}^j) - \nabla f_j(\x_t^j)\Big\rangle\right]
\n\\
\stackrel{(i)}{=}&~\frac{1}{n^2}\sum_{i=1}^n\E\left[\left\|\g_i(\x_t^i,\X_{t}^i) - \nabla f_i(\x_t^i)\right\|^2\right] \n\\
\leq&~\frac{\ol{\nu}^2}{n}.
\end{align}
We note that~$(i)$ in~\eqref{st} uses that whenever~$i\neq j$,
\begin{align}\label{ip=0}
&\E\left[\Big\langle\g_i(\x_t^i,\X_{t}^i) - \nabla f_i(\x_t^i),\g_j(\x_t^j,\X_{t}^j) - \nabla f_j(\x_t^j)\Big\rangle\big|\F_t\right] \n\\
\stackrel{(ii)}{=}&~\E\left[\Big\langle\E\left[\g_i(\x_t^i,\X_{t}^i)|\sigma(\X_t^j,\F_t)\right] - \nabla f_i(\x_t^i),\g_j(\x_t^j,\X_{t}^j) - \nabla f_j(\x_t^j)\Big\rangle\big|\F_t\right] \n\\
\stackrel{(iii)}{=}&~\E\left[\Big\langle\E\left[\g_i(\x_t^i,\X_{t}^i)|\F_t\right] - \nabla f_i(\x_t^i),\g_j(\x_t^j,\X_{t}^j) - \nabla f_j(\x_t^j)\Big\rangle\Big|\F_t\right]  \n\\
=&~0,  
\end{align}
where~$(ii)$ is due to the tower property of the conditional expectation and~$(iii)$ uses that~$\X_t^j$ is independent of~$\{\X_t^i,\F_t\}$ and~$\x_t^i$ is~$\F_t$-measurable.
Towards the third term~\eqref{vr_b_a2}, we define
for the ease of exposition,~$\forall t\geq1$,
$$\wh{\nabla}_t^i := \nabla f_i(\x_{t}^i)
- \nabla f_i(\x_{t-1}^i)$$  and recall that~$\E\left[\g_i(\x_t^i,\X_{t}^i) - \g_i(\x_{t-1}^i,\X_{t}^i)|\F_t\right] = \wh{\nabla}_t^i$.
Observe that~$\forall t\geq1$,
\begin{align}\label{zt0}
\E\left[\|\z_t\|^2|\F_t\right] 
=&~\E\bigg[\bigg\|\frac{1}{n}\sum_{i=1}^n\Big(\g_i(\x_t^i,\X_{t}^i) - \g_i(\x_{t-1}^i,\X_{t}^i) - \wh{\nabla}_t^i\Big)\bigg\|^2\Big|\F_t\bigg]
\n\\
=&~\frac{1}{n^2}\sum_{i=1}^n\E\left[\left\|\g_i(\x_t^i,\X_{t}^i) - \g_i(\x_{t-1}^i,\X_{t}^i) - \wh{\nabla}_t^i\right\|^2\big|\F_t\right] \n\\
&+\frac{1}{n^2}\sum_{i\neq j}\underbrace{\E\left[\Big\langle \g_i(\x_t^i,\X_{t}^i) - \g_i(\x_{t-1}^i,\X_{t}^i) - \wh{\nabla}_t^i, \g_j(\x_t^j,\X_{t}^j) - \g_j(\x_{t-1}^j,\X_{t}^j) - \wh{\nabla}_t^j\Big\rangle\big|\F_t\right]}_{=0}
\n\\
\stackrel{(i)}{=}&~\frac{1}{n^2}\sum_{i=1}^n\E\left[\left\|\g_i(\x_t^i,\X_{t}^i) - \g_i(\x_{t-1}^i,\X_{t}^i) - \wh{\nabla}_t\right\|^2\big|\F_t\right],  \n\\
\stackrel{(ii)}{\leq}&~\frac{1}{n^2}\sum_{i=1}^n\E\left[\left\|\g_i(\x_t^i,\X_{t}^i) - \g_i(\x_{t-1}^i,\X_{t}^i)\right\|^2\Big|\F_t\right],
\end{align}
where~$(i)$ follows from a similar line of arguments as~\eqref{ip=0} and~$(ii)$ uses the conditional variance decomposition, i.e., for any random vector~$\mb{a}\in\R^p$ consisted of square-integrable random variables,
\begin{align}\label{vrdc}
\E\left[\Big\|\mb{a}-\E\left[\mb{a}|\F_t\right]\Big\|^2|\F_t\right]
= \E\left[\left\|\mb{a}\right\|^2|\F_t\right] - \left\|\E\left[\mb{a}|\F_t\right]\right\|^2.
\end{align}
To proceed from~\eqref{zt0}, we take its expectation and observe that~$\forall t\geq1$,
\begin{align}\label{zt}
\E\left[\|\z_t\|^2\right] 
\leq&~\frac{1}{n^2}\sum_{i=1}^n\E\left[\left\|\g_i(\x_t^i,\X_{t}^i) - \g_i(\x_{t-1}^i,\X_{t}^i)\right\|^2\right]
\n\\
\stackrel{(i)}{\leq}&~\frac{L^2}{n^2}\sum_{i=1}^n\E\left[\left\|\x_t^i - \x_{t-1}^i\right\|^2\right] \n\\
=&~\frac{L^2}{n^2}\E\left[\left\|\x_t - \x_{t-1}\right\|^2\right] \n\\
=&~\frac{L^2}{n^2}\E\left[\left\|\x_t - \J\x_t
+\J\x_t - \J\x_{t-1} + \J\x_{t-1} - \x_{t-1}\right\|^2\right] \n\\
\leq&~\frac{3L^2}{n^2}\E\left[\left\|\x_t - \J\x_t\right\|^2
+n\left\|\ol{\x}_t - \ol{\x}_{t-1}\right\|^2
+\left\|\x_{t-1}- \J\x_{t-1}\right\|^2\right] \n\\
\stackrel{(ii)}{=}&~
\frac{3L^2\a^2}{n}\E\left[\left\|\ol{\mb{v}}_{t-1}\right\|^2\right]
+\frac{3L^2}{n^2}\Big(\E\left[\left\|\x_t - \J\x_t\right\|^2 + \left\|\x_{t-1}- \J\x_{t-1}\right\|^2\right]\Big),
\end{align}
where~$(i)$ uses the mean-squared smoothness of each~$\g_i$ and~$(ii)$ uses the update of~$\ol{\x}_t$ in~\eqref{ave}. The proof follows by taking the expectation~\eqref{vr_b_a2} and then using~\eqref{st} and~\eqref{zt} in the resulting inequality.

\subsection{Proof of Eq.~\eqref{vrb_sum}}
We recall from~\eqref{vr_b0} the following relationship:~$\forall t\geq1$,
\begin{align}\label{vt_recall}
\mb{v}_{t}^i - \nabla f_i(\x_t^i)
=&~\beta\Big(\g_i(\x_t^i,\X_{t}^i) - \nabla f_i(\x_t^i)\Big)
+ (1-\beta)\Big(\g_i(\x_t^i,\X_{t}^i) - \g_i(\x_{t-1}^i,\X_{t}^i) + \nabla f_i(\x_{t-1}^i)
- \nabla f_i(\x_t^i)\Big) \n\\
&+(1-\beta)\Big(\mb{v}_{t-1}^i - \nabla f_i(\x_{t-1}^i)\Big).
\end{align}
We recall that the conditional expectation of the first and the second term in~\eqref{vt_recall} with respect to~$\F_t$ is zero and that the third term in~\eqref{vt_recall} is~$\F_t$-measurable. 
Following a similar procedure in the proof of~\eqref{vr_b_a2}, we have:~$\forall t\geq1$,
\begin{align}\label{vr_b_s1}
\E\left[\|\mb{v}_{t}^i - \nabla f_i(\x_t^i)\|^2|\F_t\right]
\leq&~(1-\beta)^2\left\|\mb{v}_{t-1}^i - \nabla f_i(\x_{t-1}^i)\right\|^2 +2\beta^2\E\left[\left\|\g_i(\x_t^i,\X_{t}^i) - \nabla f_i(\x_t^i)\right\|^2|\F_t\right] \n\\
&+2(1-\beta)^2\E\left[\left\|\g_i(\x_t^i,\X_{t}^i) - \g_i(\x_{t-1}^i,\X_{t}^i) - \left( \nabla f_i(\x_{t}^i)
- \nabla f_i(\x_{t-1}^i) \right)\right\|^2|\F_t\right] \n\\   
\stackrel{(i)}{\leq}&~(1-\beta)^2\left\|\mb{v}_{t-1}^i - \nabla f_i(\x_{t-1}^i)\right\|^2 +2\beta^2\E\left[\left\|\g_i(\x_t^i,\X_{t}^i) - \nabla f_i(\x_t^i)\right\|^2|\F_t\right] \n\\
&+2(1-\beta)^2\E\left[\left\|\g_i(\x_t^i,\X_{t}^i) - \g_i(\x_{t-1}^i,\X_{t}^i) \right\|^2|\F_t\right] 
\end{align}
where~$(i)$ uses the conditional variance decomposition~\eqref{vrdc}.
We then take the expectation of~\eqref{vr_b_s1} with the help of the mean-squared smoothness and the bounded variance of each~$\g_i$ to proceed:~$\forall t\geq1$,
\begin{align}\label{vr_b_s2}
\E\left[\left\|\mb{v}_{t}^i - \nabla f_i(\x_t^i)\right\|^2\right]
\leq&~(1-\beta)^2\E\left[\left\|\mb{v}_{t-1}^i - \nabla f_i(\x_{t-1}^i)\right\|^2\right] +2\beta^2\nu_i^2 
+2(1-\beta)^2L^2\E\left[\left\|\x_t^i - \x_{t-1}^i\right\|^2\right] \n\\
\leq&~(1-\beta)^2\E\left[\left\|\mb{v}_{t-1}^i - \nabla f_i(\x_{t-1}^i)\right\|^2\right] +2\beta^2\nu_i^2 \n\\
&+6(1-\beta)^2L^2\left(\E\left[\left\|\x_t^i - \ol{\x}_{t}\right\|^2 + \left\|\ol{\x}_{t} - \ol{\x}_{t-1}\right\|^2 + \left\|\ol{\x}_{t-1} - \x_{t-1}^i\right\|^2\right]\right), \n\\
=&~(1-\beta)^2\E\left[\left\|\mb{v}_{t-1}^i - \nabla f_i(\x_{t-1}^i)\right\|^2\right] +2\beta^2\nu_i^2 
+6(1-\beta)^2L^2\a^2\E\left[\left\|\ol{\mb{v}}_{t-1}\right\|^2\right]
\n\\
&+6(1-\beta)^2L^2\E\left[\left\|\x_t^i - \ol{\x}_{t}\right\|^2 + \left\|\x_{t-1}^i-\ol{\x}_{t-1}\right\|^2\right],
\end{align}
where the last line uses the~$\ol{\x}_t$-update in~\eqref{ave}.
Summing up~\eqref{vr_b_s2} over~$i$ from~$1$ to~$n$ completes the proof.

\section{Proof of Lemma~\ref{gt_bound}}\label{proof_gt_bound}
\subsection{Proof of Lemma~\ref{gt_bound}(\ref{y1})}
Recall the initialization of~\HS~that~$\mb{v}_{-1} = \mb{0}_{np}$,~$\y_{0} = \mb{0}_{np}$, and~$\mb{v}_0^i = \frac{1}{b_0}\sum_{r=1}^{b_0}\g_i(\x_0^i,\X_{0,r}^i)$. Using the gradient tracking update~\eqref{y_up} at iteration~$t = 0$, we have:
\begin{align}\label{y10}
\E\left[\left\|\mb{y}_1-\J\y_{1}\right\|^2\right]
=&~\E\left[\left\|\W\left(\mb{y}_{0} + \mb{v}_0 - \mb{v}_{-1}\right) -\J\W\left(\mb{y}_{0} + \mb{v}_{0} - \mb{v}_{-1}\right)\right\|^2\right] \nonumber\\
\stackrel{(i)}{=}&~\E\left[\left\|\left(\W-\J\right)\mb{v}_{0}\right\|^2\right] \n\\
\stackrel{(ii)}{\leq}&~\lambda^2\E\left[\left\|\mb{v}_{0}-\nabla\mb{f}(\x_{0})+\nabla\mb{f}(\x_{0})\right\|^2\right] \n\\
=&~\lambda^2\sum_{i=1}^n\E\left[\left\|\mb{v}_{0}^i - \nabla f_i(\x_0^i)\right\|^2\right] + \lambda^2\left\|\nabla\mb{f}(\x_{0})\right\|^2 \n\\
\stackrel{(iii)}{=}&~\lambda^2\sum_{i=1}^n\E\bigg[\bigg\|\frac{1}{b_0}\sum_{r=1}^{b_0}\Big(\g_i(\x_0^i,\X_{0,r}^i) - \nabla f_i(\x_0^i)\Big)\bigg\|^2\bigg] 
+ \lambda^2\left\|\nabla\mb{f}(\x_{0})\right\|^2
\n\\
\stackrel{(iv)}{=}&~\frac{\lambda^2}{b_0^2}\sum_{i=1}^n\sum_{r=1}^{b_0}\E\left[\left\|\g_i(\x_0^i,\X_{0,r}^i) - \nabla f_i(\x_0^i)\right\|^2\right]
+\lambda^2\left\|\nabla\mb{f}(\x_{0})\right\|^2,
\end{align}
where~$(i)$ uses~$\J\W = \J$ and the initial condition of~$\mb{v}_{-1}$ and~$\mb{y}_0$,~$(ii)$ uses~$\|\W-\J\| = \lambda$,~$(iii)$ is due to the initialization of~$\mb{v}_0^i$, and~$(iv)$ follows from the fact that~$\{\X_{0,1}^i,\X_{0,2}^i,\cdots,\X_{0,b_0}^i\}$,~$\forall i\in\mc{V}$, is an independent family of random vectors, by a similar line of arguments in~\eqref{st} and~\eqref{ip=0}. The proof then follows by using the bounded variance of each~$\g_i$ in~\eqref{y10}.

\subsection{Proof of Lemma~\ref{gt_bound}(\ref{yt})}
Following the gradient tracking update~\eqref{y_up}, we have:~$\forall t\geq1$,
\begin{align}\label{yt0}
\left\|\mb{y}_{t+1}-\J\mb{y}_{t+1}\right\|^2
=&~\left\|\W\left(\mb{y}_{t} + \mb{v}_{t} - \mb{v}_{t-1}\right) -\J\W\left(\mb{y}_{t} + \mb{v}_{t} - \mb{v}_{t-1}\right)\right\|^2 \nonumber\\
\stackrel{(i)}{=}&~\left\|\W\mb{y}_{t}-\J\mb{y}_{t}+\left(\W-\J\right)\left(\mb{v}_{t} - \mb{v}_{t-1}\right)\right\|^2 \n\\
=&~\left\|\W\mb{y}_{t}-\J\mb{y}_{t}\right\|^2
+2\big\langle\W\mb{y}_{t}-\J\mb{y}_{t},\left(\W-\J\right)\left(\mb{v}_{t} - \mb{v}_{t-1}\right)\big\rangle
+\left\|\left(\W-\J\right)\left(\mb{v}_{t} - \mb{v}_{t-1}\right)\right\|^2
\n\\
\stackrel{(ii)}{\leq}&~\lambda^2\left\|\mb{y}_{t}-\J\mb{y}_{t}\right\|^2
+\underbrace{2\big\langle\W\mb{y}_{t}-\J\mb{y}_{t},\left(\W-\J\right)\left(\mb{v}_{t} - \mb{v}_{t-1}\right)\big\rangle}_{=:A_t}
+\lambda^2\left\|\mb{v}_{t} - \mb{v}_{t-1}\right\|^2,
\end{align}
where~$(i)$ uses~$\J\W=\J$ and~$(ii)$ is due to~$\|\W-\J\| = \lambda$. In the following, we bound~$A_t$ and the last term in~\eqref{yt0} respectively.
We recall the update of each local stochastic gradient estimator~$\mb{v}_t^i$ in~\eqref{vi}:~$\forall t\geq1$,
\begin{align*}
\mb{v}_t^i 
=&~\g_i(\x_t^i,\X_t^i) + (1-\beta)\mb{v}_{t-1}^i -(1-\beta)\g_i(\x_{t-1}^i,\X_t^i).
\end{align*}
We observe that~$\forall t\geq1$ and~$\forall i\in\mc{V}$,
\begin{align}\label{vt_diff}
\mb{v}_t^i - \mb{v}_{t-1}^i
=&~\g_i(\x_t^i,\X_t^i) 
-\beta\mb{v}_{t-1}^i 
-(1-\beta)\g_i(\x_{t-1}^i,\X_t^i) \n\\
=&~\g_i(\x_t^i,\X_t^i) - \g_i(\x_{t-1}^i,\X_t^i)
-\beta\mb{v}_{t-1}^i 
+\beta\g_i(\x_{t-1}^i,\X_t^i) \n\\
=&~\g_i(\x_t^i,\X_t^i) - \g_i(\x_{t-1}^i,\X_t^i)
-\beta\Big(\mb{v}_{t-1}^i 
-\nabla f_i(\x_{t-1}^i)\Big)
+\beta\Big(\g_i(\x_{t-1}^i,\X_t^i)-\nabla f_i(\x_{t-1}^i)\Big).
\end{align}
Moreover, we observe from~\eqref{vt_diff} that~$\forall t\geq1$,
\begin{align}\label{vt_diff_Ft}
\E\left[\mb{v}_t - \mb{v}_{t-1}|\F_t\right]
= \nabla\f(\x_t) - \nabla \f(\x_{t-1})
-\beta\Big(\mb{v}_{t-1} 
-\nabla\f(\x_{t-1})\Big).
\end{align}
Towards~$A_t$, we have:~$\forall t\geq1$,
\begin{align}\label{At}
\E\left[A_t|\F_t\right]    
\stackrel{(i)}{=}&~2\Big\langle\W\mb{y}_{t}-\J\mb{y}_{t},\left(\W-\J\right)\E\left[\mb{v}_{t} - \mb{v}_{t-1}|\F_t\right]\Big\rangle \n\\
\stackrel{(ii)}{=}&~2\Big\langle\W\mb{y}_{t}-\J\mb{y}_{t},\left(\W-\J\right)\Big(\nabla\f(\x_t) - \nabla \f(\x_{t-1})
-\beta\big(\mb{v}_{t-1} 
-\nabla\f(\x_{t-1})\Big)\Big\rangle \n\\
\stackrel{(iii)}{\leq}&~2\lambda\left\|\mb{y}_{t}-\J\mb{y}_{t}\right\|\cdot\lambda\left\|\nabla\f(\x_t) - \nabla \f(\x_{t-1})
-\beta\Big(\mb{v}_{t-1} 
-\nabla\f(\x_{t-1})\Big)\right\|
\n\\
\stackrel{(iv)}{\leq}&~\frac{1-\lambda^2}{2}\left\|\mb{y}_{t}-\J\mb{y}_{t}\right\|^2
+\frac{2\lambda^4}{1-\lambda^2}\left\|\nabla\f(\x_t) - \nabla \f(\x_{t-1})
-\beta\Big(\mb{v}_{t-1} 
-\nabla\f(\x_{t-1})\Big)\right\|^2,
\n\\
\stackrel{(v)}{\leq}&~\frac{1-\lambda^2}{2}\left\|\mb{y}_{t}-\J\mb{y}_{t}\right\|^2
+\frac{4\lambda^4L^2}{1-\lambda^2}\left\|\x_t - \x_{t-1}\right\|^2
+\frac{4\lambda^4\beta^2}{1-\lambda^2}\left\|\mb{v}_{t-1} 
-\nabla\f(\x_{t-1})\right\|^2,
\end{align}
where~$(i)$ is due to the~$\F_t$-measurability of~$\y_t$,~$(ii)$ uses~\eqref{vt_diff_Ft},~$(iii)$ is due to the Cauchy-Schwarz inequality and~$\|\W-\J\| = \lambda$,~$(iv)$ uses the elementary inequality that~$2ab\leq\eta a^2 + b^2/\eta$, with~$\eta=\frac{1-\lambda^2}{2\lambda^2}$ for any~$a,b\in\R$, and~$(v)$ holds since each~$f_i$ is~$L$-smooth. Next, towards the last term in~\eqref{yt0}, we take the expectation of~\eqref{vt_diff} to obtain:~$\forall t\geq1$ and~$\forall i\in\mc{V}$,
\begin{align}\label{vt_diff_norm}
\E\left[\left\|\mb{v}_t^i - \mb{v}_{t-1}^i\right\|^2\right]
\leq&~3\E\left[\left\|\g_i(\x_t^i,\X_t^i) - \g_i(\x_{t-1}^i,\X_t^i)\right\|^2\right]
+3\beta^2\E\left[\left\|\mb{v}_{t-1}^i 
-\nabla f_i(\x_{t-1}^i)\right\|^2\right] \n\\
&+3\beta^2\E\left[\left\|\g_i(\x_{t-1}^i,\X_t^i)-\nabla f_i(\x_{t-1}^i)\right\|^2\right]\n\\
\leq&~3L^2\E\left[\left\|\x_t^i-\x_{t-1}^i\right\|^2\right] +3\beta^2\E\left[\left\|\mb{v}_{t-1}^i 
-\nabla f_i(\x_{t-1}^i)\right\|^2\right] + 3\beta^2\nu_i^2,
\end{align}
where~\eqref{vt_diff_norm} is due to the mean-squared smoothness and the bounded variance of each~$\g_i$. Summing up~\eqref{vt_diff_norm} over~$i$ from~$1$ to~$n$ gives an upper bound on the last term in~\eqref{yt0}:~$\forall t\geq1$,
\begin{align}\label{vt_diff_norm_sum}
\lambda^2\E\left[\left\|\mb{v}_t - \mb{v}_{t-1}\right\|^2\right]
\leq&~
3\lambda^2L^2\E\left[\left\|\x_t - \x_{t-1}\right\|^2\right] +3\lambda^2\beta^2\E\left[\left\|\mb{v}_{t-1} 
-\nabla\f(\x_{t-1})\right\|^2\right] + 3\lambda^2n\beta^2\ol{\nu}^2. 
\end{align}
We now use~\eqref{At} and~\eqref{vt_diff_norm_sum} in~\eqref{yt0} to obtain:~$\forall t\geq1$,
\begin{align}\label{yt1}
\E\left[\left\|\mb{y}_{t+1}-\J\mb{y}_{t+1}\right\|^2\right] 
\leq&~
\frac{1+\lambda^2}{2}\E\left[\left\|\mb{y}_{t}-\J\mb{y}_{t}\right\|^2\right] 
+ \frac{7\lambda^2L^2}{1-\lambda^2}\E\left[\left\|\x_t-\x_{t-1}\right\|^2\right]
\n\\
&
+\frac{7\lambda^2\beta^2}{1-\lambda^2}\E\left[\left\|\mb{v}_{t-1}
-\nabla\f(\x_{t-1})\right\|^2\right] 
+ 3\lambda^2n\beta^2\ol{\nu}^2.
\end{align}
Towards the second term in~\eqref{yt1}, we use~\eqref{consensus_2} to obtain:~$\forall t\geq1$,
\begin{align}\label{xt_diff_2}
\|\x_t-\x_{t-1}\|^2
=&~\left\|\x_t-\J\x_t + \J\x_{t} - \J\x_{t-1} + \J\x_{t-1} -\x_{t-1}\right\|^2 \n\\
\stackrel{(i)}{\leq}&~3\left\|\x_t-\J\x_t\right\|^2 + 3n\a^2\left\|\ol{\mb{v}}_{t-1}\right\|^2
+ 3\left\|\x_{t-1}-\J\x_{t-1}\right\|^2 \n\\
\leq&~6\lambda^2\a^2\left\|\y_t-\J\y_t\right\|^2 + 3n\a^2\left\|\ol{\mb{v}}_{t-1}\right\|^2
+ 9\left\|\x_{t-1}-\J\x_{t-1}\right\|^2,
\end{align}
where~$(i)$ uses the~$\ol{\x}_t$-update in~\eqref{ave}.
Finally, we use~\eqref{xt_diff_2} in~\eqref{yt1} to obtain:~$\forall t\geq1$,
\begin{align*}
\E\left[\left\|\mb{y}_{t+1}-\J\mb{y}_{t+1}\right\|^2\right] 
\leq&~
\left(\frac{1+\lambda^2}{2} + \frac{42\lambda^4L^2\a^2}{1-\lambda^2}\right)\E\left[\left\|\mb{y}_{t}-\J\mb{y}_{t}\right\|^2\right] 
+ \frac{21\lambda^2nL^2\a^2}{1-\lambda^2}\E\left[\left\|\ol{\mb{v}}_{t-1}\right\|^2\right]
\n\\
&
+ \frac{63\lambda^2L^2}{1-\lambda^2}\E\left[\left\|\x_{t-1}-\J\x_{t-1}\right\|^2\right]
+\frac{7\lambda^2\beta^2}{1-\lambda^2}\E\left[\left\|\mb{v}_{t-1}
-\nabla\f(\x_{t-1})\right\|^2\right] 
+ 3\lambda^2n\beta^2\ol{\nu}^2.
\end{align*}
The proof is completed by the fact that~$\frac{1+\lambda^2}{2}+\frac{42\lambda^4L^2\a^2}{1-\lambda^2}\leq\frac{3+\lambda^2}{4}$ if~$0<\a\leq\frac{1-\lambda^2}{2\sqrt{42}\lambda^2L}$.  

\section{Proof of Lemma~\ref{convs}}\label{proof_convs}
\subsection{Proof of Eq.~\eqref{convs1}}
We recursively apply the inequality on~$V_t$ from~$t$ to~$0$ to obtain:~$\forall t\geq1$,
\begin{align}\label{recur_base_0}
V_t
\leq&~qV_{t-1} + qR_{t-1} + Q_t + C \n\\
\leq&~q^2V_{t-2} + (q^2R_{t-2} + qR_{t-1}) 
+ (qQ_{t-1} + Q_t)
+ (qC+ C) \n\\ 
\cdots& \n\\
\leq&~q^tV_0 
+ \sum_{i = 0}^{t-1}q^{t-i}R_{i} 
+ \sum_{i = 1}^{t}q^{t-i}Q_{i} 
+ C\sum_{i=0}^{t-1}q^i.
\end{align}
Summing up~\eqref{recur_base_0} over~$t$ from~$1$ to~$T$ gives:~$\forall T\geq1$,
\begin{align*}
\sum_{t=0}^{T}V_t
\leq&~V_0\sum_{t=0}^{T}q^t 
+ \sum_{t=1}^{T}\sum_{i = 0}^{t-1}q^{t-i}R_{i} 
+ \sum_{t=1}^{T}\sum_{i = 1}^{t}q^{t-i}Q_{i}
+ C\sum_{t=1}^{T}\sum_{i=0}^{t-1}q^i \n\\
\leq&~V_0\sum_{t=0}^{\infty}q^t 
+ \sum_{t=0}^{T-1}\left(\sum_{i = 0}^{\infty}q^{i}\right)R_{t} 
+ \sum_{t=1}^{T}\left(\sum_{i = 0}^{\infty}q^{i}\right)Q_{t}
+ C\sum_{t=1}^{T}\sum_{i=0}^{\infty}q^i,
\end{align*}
and the proof follows by~$\sum_{i=0}^{\infty}q^i = (1-q)^{-1}$.

\subsection{Proof of Eq.~\eqref{convs2}}
We recursively apply the inequality on~$V_t$ from~$t+1$ to~$1$ to obtain:~$\forall t\geq1$,
\begin{align}\label{recur_base_tt0}
V_{t+1} 
\leq&~qV_{t} + R_{t-1} + C         \n\\
\leq&~q^2V_{t-1} + (qR_{t-2} + R_{t-1}) + (qC + C) \n\\
\cdots& \n\\
\leq&~q^{t}V_{1} + \sum_{i=0}^{t-1}q^{t-1-i}R_i
+C\sum_{i=0}^{t-1}q^i.
\end{align}
We sum up~\eqref{recur_base_tt0} over~$t$ from~$1$ to~$T-1$ to obtain:~$\forall T\geq2$,
\begin{align*}
\sum_{t=0}^{T-1}V_{t+1} 
\leq&~V_{1}\sum_{t=0}^{T-1}q^{t} + \sum_{t=1}^{T-1}\sum_{i=0}^{t-1}q^{t-1-i}R_i
+C\sum_{t=1}^{T-1}\sum_{i=0}^{t-1}q^i       \n\\   
\leq&~V_{1}\sum_{t=0}^{\infty}q^{t} + \sum_{t=0}^{T-2}\left(\sum_{i=0}^{\infty}q^{i}\right)R_t
+C\sum_{t=1}^{T-1}\sum_{i=0}^{\infty}q^i,
\end{align*}
and the proof follows by~$\sum_{i=0}^{\infty}q^i = (1-q)^{-1}$.

\section{Proof of Lemma~\ref{vrb_acc}}\label{proof_vrb_acc}
\subsection{Proof of Eq.~\eqref{vrb_a_acc}}
We first observe that~$\frac{1}{1-(1-\beta)^2} \leq\frac{1}{\beta}$ for~$\beta\in(0,1)$. Applying~\eqref{convs1} to~\eqref{vrb_a} gives:~$\forall T\geq1$,
\begin{align}\label{vrb_a_acc_0}
&\sum_{t=0}^{T}\E\left[\left\|\ol{\mb{v}}_{t} - \ol{\nabla \f}(\x_t)\right\|^2\right]  \n\\
\leq&~\frac{\E\left[\|\ol{\mb{v}}_{0} - \ol{\nabla \f}(\x_{0})\|^2\right]}{\beta}
+\frac{6L^2\a^2}{n\beta}\sum_{t=0}^{T-1}\E\left[\left\|\ol{\mb{v}}_{t}\right\|^2\right] 
+\frac{6L^2}{n^2\beta}\sum_{t=0}^{T-1}\E\left[\left\|\x_{t+1} - \J\x_{t+1}\|^2 + \|\x_{t}- \J\x_{t}\right\|^2\right] +\frac{2\beta\ol{\nu}^2 T}{n} \n\\
\leq&~\frac{\E\left[\|\ol{\mb{v}}_{0} - \ol{\nabla \f}(\x_{0})\|^2\right]}{\beta}
+\frac{6L^2\a^2}{n\beta}\sum_{t=0}^{T-1}\E\left[\left\|\ol{\mb{v}}_{t}\right\|^2\right] 
+\frac{12L^2}{n^2\beta}\sum_{t=0}^{T}\E\left[\left\|\x_{t}- \J\x_{t}\right\|^2\right] +\frac{2\beta\ol{\nu}^2 T}{n}.
\end{align}
Towards the first term in~\eqref{vrb_a_acc_0}, we observe that
\begin{align}\label{v0_b}
\E\left[\left\|\ol{\mb{v}}_0-\ol{\nabla\f}(\x_0)\right\|^2\right]
=&~\E\left[\Big\|\frac{1}{n}\sum_{i=1}^n\frac{1}{b_0}\sum_{r=1}^{b_0}\Big(\g_i(\x_0^i,\X_{0,r}^i) -\nabla f_i(\x_0^i)\Big)\Big\|^2\right] \n\\
\stackrel{(i)}{=}&~\frac{1}{n^2b_0^2}\sum_{i=1}^n\sum_{r=1}^{b_0}\E\left[\left\|\g_i(\x_0^i,\X_{0,r}^i) -\nabla f_i(\x_0^i)\right\|^2\right] \leq\frac{\ol{\nu}^2}{nb_0},    
\end{align}
where~$(i)$ follows from a similar line of arguments in~\eqref{ip=0}. Then~\eqref{vrb_a_acc} follows from using~\eqref{v0_b} in~\eqref{vrb_a_acc_0}.

\subsection{Proof of Eq.~\eqref{vrb_s_acc}}
We apply~\eqref{convs1} to~\eqref{vrb_sum} to obtain:~$\forall T\geq1$,
\begin{align}\label{vrb_s_acc_0}
&\sum_{t=0}^{T}\E\left[\left\|\mb{v}_{t} - \nabla\f(\x_t)\right\|^2\right]  \n\\
\leq&~\frac{\E\left[\|\mb{v}_{0} - \nabla \f(\x_{0})\|^2\right]}{\beta}
+\frac{6nL^2\a^2}{\beta}\sum_{t=0}^{T-1}\E\left[\left\|\ol{\mb{v}}_{t}\right\|^2\right] 
+\frac{6L^2}{\beta}\sum_{t=0}^{T-1}\E\left[\left\|\x_{t+1} - \J\x_{t+1}\|^2 + \|\x_{t}- \J\x_{t}\right\|^2\right] +2n\beta\ol{\nu}^2T \n\\
\leq&~\frac{\E\left[\|\mb{v}_{0} - \nabla \f(\x_{0})\|^2\right]}{\beta}
+\frac{6nL^2\a^2}{\beta}\sum_{t=0}^{T-1}\E\left[\left\|\ol{\mb{v}}_{t}\right\|^2\right] 
+\frac{12L^2}{\beta}\sum_{t=0}^{T}\E\left[\left\|\x_{t}- \J\x_{t}\right\|^2\right] + 2n\beta\ol{\nu}^2 T.
\end{align}
In~\eqref{vrb_s_acc_0}, we observe that
\begin{align}\label{v0_b_s}
\E\left[\|\mb{v}_0-\nabla\f(\x_0)\|^2\right]
=&~\sum_{i=1}^n\E\left[\Big\|\frac{1}{b_0}\sum_{r=1}^{b_0}\left(\g_i(\x_0^i,\X_{0,r}^i) -\nabla f_i(\x_0^i)\right)\Big\|^2\right] \n\\
\stackrel{(i)}{=}&~\frac{1}{b_0^2}\sum_{i=1}^n\sum_{r=1}^{b_0}\E\left[\left\|\g_i(\x_0^i,\X_{0,r}^i) -\nabla f_i(\x_0^i)\right\|^2\right] \leq \frac{n\ol{\nu}^2}{b_0},    
\end{align}
where~$(i)$ follows from a similar line of arguments in~\eqref{ip=0}. Then~\eqref{vrb_s_acc} follows from using~\eqref{v0_b_s} in~\eqref{vrb_s_acc_0}.

\section{Proof of Lemma~\ref{cons_bound_sum_final}}\label{proof_cons}
We recall that~$\|\x_t-\J\x_t\| = 0$, since it is assumed without generality that~$\x_0^i = \x_0^j$ for any~$i,j\in\mc{V}$. Applying~\eqref{convs1} to~\eqref{consensus_1} yields:~$\forall T\geq1$,
\begin{align}\label{y->x}
\sum_{t=0}^{T}\left\|\mb{x}_{t}-\J\mb{x}_{t}
\right\|^2 
\leq  
\frac{4\lambda^2\alpha^2}{(1-\lambda^2)^2}\sum_{t=1}^{T}\left\|\mb{y}_{t}-\J\mb{y}_{t}\right\|^2.    
\end{align}
To further bound~$\sum_{t=1}^{T}\left\|\mb{y}_{t}-\J\mb{y}_{t}\right\|^2$, we apply~\eqref{convs2} in Lemma~\ref{gt_bound}(\ref{yt}) to obtain: if~$0<\a\leq\frac{1-\lambda^2}{2\sqrt{42}\lambda^2L}$, then~$\forall T\geq2$,
\begin{align}\label{gt_bound_sum_0}
&\sum_{t=1}^{T}\E\left[\left\|\mb{y}_{t}-\J\mb{y}_{t}\right\|^2\right]  \n\\
\leq&~
\frac{4\E\left[\|\mb{y}_{1}-\J\mb{y}_{1}\|^2\right]}{1-\lambda^2}
+ \frac{84\lambda^2nL^2\a^2}{(1-\lambda^2)^2}\sum_{t=0}^{T-2}\E\left[\|\ol{\mb{v}}_{t}\|^2\right]
+ \frac{252\lambda^2L^2}{(1-\lambda^2)^2}\sum_{t=0}^{T-2}\E\left[\|\x_{t} - \J\x_{t}\|^2\right]
\n\\
&
+\frac{28\lambda^2\beta^2}{(1-\lambda^2)^2}\sum_{t=0}^{T-2}\E\left[\|\mb{v}_t 
-\nabla\f(\x_{t})\|^2\right]
+ \frac{12\lambda^2n\beta^2\ol{\nu}^2T}{1-\lambda^2} \n\\
\leq&~
\frac{84\lambda^2nL^2\a^2}{(1-\lambda^2)^2}\sum_{t=0}^{T-2}\E\left[\|\ol{\mb{v}}_{t}\|^2\right]
+ \frac{252\lambda^2L^2}{(1-\lambda^2)^2}\sum_{t=0}^{T-2}\E\left[\|\x_{t} - \J\x_{t}\|^2\right]
\n\\
&
+\frac{28\lambda^2\beta^2}{(1-\lambda^2)^2}\sum_{t=0}^{T-2}\E\left[\|\mb{v}_t 
-\nabla\f(\x_{t})\|^2\right]
+\frac{12\lambda^2n\beta^2\ol{\nu}^2T}{1-\lambda^2}
+\frac{4\lambda^2\left\|\nabla\mb{f}\big(\x_{0}\big)\right\|^2}{1-\lambda^2}
+ \frac{4\lambda^2n\ol{\nu}^2}{(1-\lambda^2)b_0},
\end{align}
where the last inequality is due to Lemma~\ref{gt_bound}(\ref{y1}).
To proceed, we use~\eqref{vrb_s_acc}, an upper bound on~$\sum_{t}\E\left[\|\mb{v}_{t} - \nabla\f(\x_t)\|^2\right]$, in~\eqref{gt_bound_sum_0} to obtain: if~$0<\a\leq\frac{1-\lambda^2}{2\sqrt{42}\lambda^2L}$ and~$\beta\in(0,1)$, then~$\forall T\geq2$,
\begin{align}\label{gt_bound_sum_1}
\sum_{t=1}^{T}\E\left[\left\|\mb{y}_{t}-\J\mb{y}_{t}\right\|^2\right]  
\leq&~
\frac{252\lambda^2nL^2\a^2}{(1-\lambda^2)^2}\sum_{t=0}^{T-2}\E\left[\|\ol{\mb{v}}_{t}\|^2\right]
+ \frac{588\lambda^2L^2}{(1-\lambda^2)^2}\sum_{t=0}^{T-1}\E\left[\|\x_{t} - \J\x_{t}\|^2\right]
\n\\
&
+\frac{28\lambda^2n\beta\ol{\nu}^2}{(1-\lambda^2)^2b_0}
+\frac{56\lambda^2n\beta^3\ol{\nu}^2T}{(1-\lambda^2)^2}
+\frac{12\lambda^2n\beta^2\ol{\nu}^2T}{1-\lambda^2}
+\frac{4\lambda^2\left\|\nabla\mb{f}\big(\x_{0}\big)\right\|^2}{1-\lambda^2}
+ \frac{4\lambda^2n\ol{\nu}^2}{(1-\lambda^2)b_0} \n\\
=&~
\frac{252\lambda^2nL^2\a^2}{(1-\lambda^2)^2}\sum_{t=0}^{T-2}\E\left[\|\ol{\mb{v}}_{t}\|^2\right]
+ \frac{588\lambda^2L^2}{(1-\lambda^2)^2}\sum_{t=0}^{T-1}\E\left[\|\x_{t} - \J\x_{t}\|^2\right]
\n\\
&
+ \left(\frac{7\beta}{1-\lambda^2}
+ 1\right)\frac{4\lambda^2n\ol{\nu}^2}{(1-\lambda^2)b_0}
+\left(\frac{14\beta}{1-\lambda^2}
+ 3\right)\frac{4\lambda^2n\beta^2\ol{\nu}^2T}{1-\lambda^2}
+ \frac{4\lambda^2\left\|\nabla\mb{f}\big(\x_{0}\big)\right\|^2}{1-\lambda^2}.
\end{align}
Finally, we use~\eqref{gt_bound_sum_1} in~\eqref{y->x} to obtain:~$\forall T\geq2$,
\begin{align*}
\sum_{t=0}^{T}\E\left[\left\|\mb{x}_{t}-\J\mb{x}_{t}\right\|^2\right]  
\leq&~
\frac{1008\lambda^4nL^2\a^4}{(1-\lambda^2)^4}\sum_{t=0}^{T-2}\E\left[\|\ol{\mb{v}}_{t}\|^2\right]
+ \frac{2352\lambda^4L^2\a^2}{(1-\lambda^2)^4}\sum_{t=0}^{T-1}\E\left[\|\x_{t} - \J\x_{t}\|^2\right]
\n\\
&
+ \left(\frac{7\beta}{1-\lambda^2}
+ 1\right)\frac{16\lambda^4n\ol{\nu}^2\a^2}{(1-\lambda^2)^{3}b_0}
+\left(\frac{14\beta}{1-\lambda^2}
+ 3\right)\frac{16\lambda^4n\beta^2\ol{\nu}^2\a^2T}{(1-\lambda^2)^3}
+ \frac{16\lambda^4\left\|\nabla\mb{f}\big(\x_{0}\big)\right\|^2\a^2}{(1-\lambda^2)^3},
\end{align*}
which may be written equivalently as
\begin{align}\label{cons_bound_sum}
\left(1-\frac{2352\lambda^4L^2\a^2}{(1-\lambda^2)^4}\right)\sum_{t=0}^{T}\E\left[\left\|\mb{x}_{t}-\J\mb{x}_{t}\right\|^2\right]  
\leq&~
\frac{1008\lambda^4nL^2\a^4}{(1-\lambda^2)^4}\sum_{t=0}^{T-2}\E\left[\|\ol{\mb{v}}_{t}\|^2\right]
+ \left(\frac{7\beta}{1-\lambda^2}
+ 1\right)\frac{16\lambda^4n\ol{\nu}^2\a^2}{(1-\lambda^2)^{3}b_0}
\n\\
&
+\left(\frac{14\beta}{1-\lambda^2}
+ 3\right)\frac{16\lambda^4n\beta^2\ol{\nu}^2\a^2T}{(1-\lambda^2)^3}
+ \frac{16\lambda^4\left\|\nabla\mb{f}\big(\x_{0}\big)\right\|^2\a^2}{(1-\lambda^2)^3}.
\end{align}
We observe in~\eqref{cons_bound_sum} that~$\frac{2352\lambda^4L^2\a^2}{(1-\lambda^2)^4}\leq\frac{1}{2}$ if~$0<\a\leq\frac{(1-\lambda^2)^2}{70\lambda^2L}$, and the proof follows.

\section{Proof of Theorem~\ref{main0}}\label{proof_main0}
For the ease of presentation, we denote~$\Delta_0:= F(\ol{\x}_0)- F^*$ in the following. We apply~\eqref{vrb_a_acc} to Lemma~\ref{DS0} to obtain: if~$0<\a\leq\frac{1}{2L}$, then~$\forall T\geq1$,
\begin{align}\label{ds1_0}
\sum_{t=0}^{T}\E\left[\left\|\nabla F(\ol{\x}_{t})\right\|^2\right]
\leq&~ 
\frac{2\Delta_0}{\a}
- \frac{1}{2}\sum_{t=0}^{T}\E\left[\left\|\ol{\mb{v}}_{t}\right\|^2\right] 
+ \frac{2L^2}{n}\sum_{t=0}^{T}\E\left[\left\|\mb{x}_{t}-\J\x_{t}\right\|^2\right] \n\\
&+ \frac{2\ol{\nu}^2}{\beta b_0n}
+\frac{12L^2\a^2}{n\beta}\sum_{t=0}^{T-1}\E\left[\left\|\ol{\mb{v}}_{t}\right\|^2\right] 
+\frac{24L^2}{n^2\beta}
\sum_{t=0}^{T}\E\left[\left\|\x_{t} - \J\x_{t}\right\|^2\right]
+\frac{4\beta\ol{\nu}^2T}{n} 
\n\\
\leq&~ 
\frac{2\Delta_0}{\a}
- \frac{1}{4}\sum_{t=0}^{T}\E\left[\left\|\ol{\mb{v}}_{t}\right\|^2\right] 
+ \frac{2L^2}{n}\left(1+\frac{12}{n\beta}\right)\sum_{t=0}^{T}\E\left[\left\|\mb{x}_{t}-\J\x_{t}\right\|^2\right] \n\\
&+ \frac{2\ol{\nu}^2}{\beta b_0n}
+\frac{4\beta\ol{\nu}^2T}{n}
-\left(\frac{1}{4}-\frac{12L^2\a^2}{n\beta}\right)\sum_{t=0}^{T}\E\left[\left\|\ol{\mb{v}}_{t}\right\|^2\right].
\end{align}
Therefore, if~$0<\a<\frac{1}{4\sqrt{3}L}$ and~$\frac{48L^2\a^2}{n}\leq\beta<1$, i.e.,~$\frac{1}{4} - \frac{12L^2\a^2}{n\beta}\geq0$, we may drop the last term in~\eqref{ds1_0} to obtain:~$\forall T\geq1$,
\begin{align}\label{DS1}
\sum_{t=0}^{T}\E\left[\left\|\nabla F(\ol{\x}_{t})\right\|^2\right] 
\leq&~ 
\frac{2\Delta_0}{\a}
- \frac{1}{4}\sum_{t=0}^{T}\E\left[\left\|\ol{\mb{v}}_{t}\right\|^2\right] 
+ \frac{2L^2}{n}\left(1+\frac{12}{n\beta}\right)\sum_{t=0}^{T}\E\left[\left\|\mb{x}_{t}-\J\x_{t}\right\|^2\right] 
+ \frac{2\ol{\nu}^2}{\beta b_0n}
+\frac{4\beta\ol{\nu}^2T}{n}.
\end{align}
Moreover, we observe:~$\forall T\geq1$,
\begin{align}\label{local0}
\frac{1}{n}\sum_{i=1}^n\sum_{t=0}^{T}\E\left[\left\|\nabla F(\x_t^i)\right\|^2\right]    
\leq&~\frac{2}{n}\sum_{i=1}^n\sum_{t=0}^{T}\E\left[\left\|\nabla F(\x_t^i) - \nabla F(\ol{\x}_t)\right\|^2 + \left\|\nabla F(\ol{\x}_t)\right\|^2\right]\n\\
=&~\frac{2L^2}{n}\sum_{t=0}^{T}\E\left[\left\|\x_t - \J\x_t\right\|^2\right]
+2\sum_{t=0}^{T}\E\left[\left\|\nabla F(\ol{\x}_t)\right\|^2\right],
\end{align}
where the last line uses the~$L$-smoothness of~$F$. Using~\eqref{DS1} in~\eqref{local0} yields: if~$0<\a<\frac{1}{4\sqrt{3}L}$ and~$48L^2\a^2/n\leq\beta<1$, then
\begin{align}\label{local1}
\frac{1}{n}\sum_{i=1}^n\sum_{t=0}^{T}\E\left[\left\|\nabla F(\x_t^i)\right\|^2\right]    
\leq&~\frac{4\Delta_0}{\a}
- \frac{1}{2}\sum_{t=0}^{T}\E\left[\left\|\ol{\mb{v}}_{t}\right\|^2\right] 
+ \frac{6L^2}{n}\left(1+\frac{8}{n\beta}\right)\sum_{t=0}^{T}\E\left[\left\|\mb{x}_{t}-\J\x_{t}\right\|^2\right] 
+ \frac{4\ol{\nu}^2}{\beta b_0n}
+\frac{8\beta\ol{\nu}^2T}{n}.
\end{align}
According to~\eqref{local1}, if~$0<\alpha<\frac{1}{4\sqrt{3}L}$ and~$\beta = 48L^2\a^2/n$, we have:~$\forall T\geq1$,
\begin{align}\label{ds2_0}
\frac{1}{n}\sum_{i=1}^n\sum_{t=0}^{T}\E\left[\left\|\nabla F(\x_t^i)\right\|^2\right]    
\leq&~ 
\frac{4\Delta_0}{\a}
- \frac{1}{2}\sum_{t=0}^{T}\E\left[\left\|\ol{\mb{v}}_{t}\right\|^2\right] 
+ \frac{6L^2}{n}\left(1+\frac{1}{6L^2\a^2}\right)\sum_{t=0}^{T}\E\left[\left\|\mb{x}_{t}-\J\x_{t}\right\|^2\right] 
+ \frac{4\ol{\nu}^2}{\beta b_0n}
+\frac{8\beta\ol{\nu}^2T}{n} \n\\
\leq&~ 
\frac{4\Delta_0}{\a}
\underbrace{-\frac{1}{2}\sum_{t=0}^{T}\E\left[\left\|\ol{\mb{v}}_{t}\right\|^2\right] 
+ \frac{2}{n\a^2}\sum_{t=0}^{T}\E\left[\left\|\mb{x}_{t}-\J\x_{t}\right\|^2\right]}_{=:\Phi_T}
+ \frac{4\ol{\nu}^2}{\beta b_0n}
+\frac{8\beta\ol{\nu}^2T}{n},
\end{align}
where the last line is due to~$6L^2\a^2<1/8$. To simplify~$\Phi_T$, we use Lemma~\ref{cons_bound_sum_final} to obtain: if~$0<\a\leq\frac{(1-\lambda^2)^2}{70\lambda^2L}$ then~$\forall T\geq2$,
\begin{align}\label{Ut0}
\Phi_T
\leq&~
-\frac{1}{2}\left(1- 
\frac{8064\lambda^4L^2\a^2}{(1-\lambda^2)^4}\right)\sum_{t=0}^{T}\E\left[\left\|\ol{\mb{v}}_{t}\right\|^2\right]
+ \frac{64\lambda^4}{(1-\lambda^2)^3}\frac{\left\|\nabla\mb{f}\big(\x_{0}\big)\right\|^2}{n}
\n\\
&
+ \left(\frac{7\beta}{1-\lambda^2}
+ 1\right)\frac{64\lambda^4\ol{\nu}^2}{(1-\lambda^2)^{3}b_0}
+\left(\frac{14\beta}{1-\lambda^2}
+ 3\right)\frac{64\lambda^4\beta^2\ol{\nu}^2T}{(1-\lambda^2)^3}.
\end{align}
In~\eqref{Ut0}, we observe that if~$0<\a\leq\frac{(1-\lambda^2)^2}{90\lambda^2L}$, then~$1- 
\frac{8064\lambda^4L^2\a^2}{(1-\lambda^2)^4}\geq0$ and thus the first term in~\eqref{Ut0} may be dropped; moreover, if~$0<\a\leq\frac{\sqrt{n(1-\lambda^2)}}{26\lambda L}$, then $\beta = \frac{48L^2\a^2}{n} \leq \frac{1-\lambda^2}{14\lambda^2}$. Hence, if~$0<\a\leq\min\Big\{\frac{(1-\lambda^2)^2}{90\lambda^2},\frac{\sqrt{n(1-\lambda^2)}}{26\lambda}\Big\}\frac{1}{L}$, then~\eqref{Ut0} reduces to:~$\forall T\geq2$,
\begin{align}\label{Ut1}
\Phi_T
\leq&~
\frac{64\lambda^4}{(1-\lambda^2)^3}\frac{\left\|\nabla\mb{f}\big(\x_{0}\big)\right\|^2}{n}
+ \frac{96\lambda^2\ol{\nu}^2}{(1-\lambda^2)^{3}b_0}
+\frac{256\lambda^2\beta^2\ol{\nu}^2T}{(1-\lambda^2)^3}.
\end{align}
Finally, we use~\eqref{Ut1} in~\eqref{ds2_0} to obtain:
if~$0<\a<\min\Big\{\frac{1}{4\sqrt{3}},\frac{(1-\lambda^2)^2}{90\lambda^2},\frac{\sqrt{n(1-\lambda^2)}}{26\lambda}\Big\}\frac{1}{L}$, we have:~$\forall T\geq2$,
\begin{align}\label{final}
\frac{1}{n(T+1)}\sum_{i=1}^n\sum_{t=0}^{T}\E\left[\left\|\nabla F(\x_t^i)\right\|^2\right]    
\leq&~ 
\frac{4\Delta_0}{\a T}
+ \frac{4\ol{\nu}^2}{\beta b_0nT}
+\frac{8\beta\ol{\nu}^2}{n} \n\\
&+\frac{64\lambda^4}{(1-\lambda^2)^3T}\frac{\left\|\nabla\mb{f}\big(\x_{0}\big)\right\|^2}{n}
+ \frac{96\lambda^2\ol{\nu}^2}{(1-\lambda^2)^{3}b_0T}
+\frac{256\lambda^2\beta^2\ol{\nu}^2}{(1-\lambda^2)^3}.
\end{align}
The proof follows by~\eqref{final} and that~$\E[\|\nabla F(\wt{\x}_T)\|^2] = \frac{1}{n(T+1)}\sum_{i=1}^n\sum_{t=0}^{T}\E[\|\nabla F(\x_t^i)\|^2]$ since~$\wt{\x}_T$ is chosen uniformly at random from~$\{\x_t^i:\forall i\in\mc{V}, 0\leq t\leq T\}$.

\end{document}